\input amstex
\input amsppt.sty
\magnification=\magstep1
\hsize=34truecc
 \vsize=22.5truecm
\baselineskip=14truept
\NoBlackBoxes
\def\q{\quad}
\def\qq{\qquad}
\def\mod#1{\ (\text{\rm mod}\ #1)}
\def\t{\text}
\def\qtq#1{\q\t{#1}\q}

\def\f{\frac}
\def\e{\equiv}
\def\a{\alpha}

\def\ep{\varepsilon}

\def\sls#1#2{(\f{#1}{#2})}
 \def\ls#1#2{\big(\f{#1}{#2}\big)}
\def\Ls#1#2{\Big(\f{#1}{#2}\Big)}
\def\sqs#1#2{(\f{#1}{#2})_4}
\def\qs#1#2{\big(\f{#1}{#2}\big)_4}
\def\Qs#1#2{\Big(\f{#1}{#2}\Big)_4}

\let \pro=\proclaim
\let \endpro=\endproclaim
\topmatter
\title\nofrills Quartic, octic residues and binary quadratic forms
\endtitle
\author ZHI-Hong Sun\endauthor
\affil School of Mathematical Sciences, Huaiyin Normal University,
\\ Huaian, Jiangsu 223001, PR China
\\ E-mail: zhihongsun$\@$yahoo.com
\\ Homepage: http://www.hytc.edu.cn/xsjl/szh
\endaffil
 \nologo \NoRunningHeads
 \abstract{Let $\Bbb Z$ be the set of integers, and let $(m,n)$ be
 the greatest common divisor of integers $m$ and $n$.
 Let $p\equiv 1\mod
4$ be a prime, $q\in\Bbb Z$, $2\nmid q$ and $p=c^2+d^2=x^2+qy^2$
with $c,d,x,y\in\Bbb Z$ and $c\e 1\mod 4$. Suppose that $(c,x+d)=1$
or $(d,x+c)$ is a power of $2$. In the paper, by using the quartic
reciprocity law we determine $q^{[p/8]}\mod p$ in terms of $c,d,x$
and $y$, where $[\cdot]$ is the greatest integer function. We also
determine
$\big(\frac{b+\sqrt{b^2+4^{\alpha}}}2\big)^{\frac{p-1}4}\mod p$ for
odd $b$ and $(2a+\sqrt{4a^2+1})^{\f{p-1}4}\mod p$ for $a\in\Bbb Z$.
 As applications we obtain the congruence for $U_{\f{p-1}4}\mod p$
  and the criterion for $p\mid U_{\frac{p-1}8}$ (if $p\equiv 1\mod 8$),
   where $\{U_n\}$ is the
Lucas sequence given by $U_0=0,\ U_1=1$ and $U_{n+1}=bU_n+U_{n-1}\
(n\ge 1)$, and $b\not\equiv 2\mod 4$. Hence we partially solve some
conjectures posed by the author in two previous papers.
\par\q
}
 \endabstract
 \footnote"" {The author is
supported by the National Natural Sciences Foundation of China (No.
10971078).}
\endtopmatter
\document
\subheading{1. Introduction}
\par Let $\Bbb Z$ and $\Bbb
N$ be the set of integers and the set of positive integers,
respectively, $i=\sqrt{-1}$ and $\Bbb Z[i]=\{a+bi\mid a,b\in\Bbb
Z\}$. For any positive odd number $m$ and $a\in\Bbb Z$ let $(\f am)$
be the (quadratic) Jacobi symbol. (We also assume $(\f a1)=1$.) For
our convenience we also define $(\f a{-m})=(\f am)$. Then for any
two odd numbers $m$ and $n$ with $m>0$ or $n>0$ we have the
following general quadratic reciprocity law: $\sls
mn=(-1)^{\f{m-1}2\cdot \f{n-1}2}\sls nm$.
\par  For $a,b,c,d\in\Bbb Z$ with $2\nmid c$ and $2\mid d$, one can define the
quartic Jacobi symbol $\qs{a+bi}{c+di}$ as in [S5]. From [IR] we
know that
$\overline{\qs{a+bi}{c+di}}=\qs{a-bi}{c-di}=\qs{a+bi}{c+di}^{-1}$,
where $\bar x$ means the complex conjugate of $x$.
 In Section 2 we list main properties of the quartic Jacobi
symbol. See also [IR], [BEW] and [S3].
\par For $b,c\in\Bbb Z$ the Lucas sequences
$\{U_n(b,c)\}$ and $\{V_n(b,c)\}$ are defined by
$$\aligned &U_0(b,c)=0,\ U_1(b,c)=1,\\& U_{n+1}(b,c)=bU_n(b,c)-cU_{n-1}(b,c)\
(n\ge 1)\endaligned$$ and $$\aligned &V_0(b,c)=2,\ V_1(b,c)=b,\\&
V_{n+1}(b,c)=bV_n(b,c)-cV_{n-1}(b,c)\ (n\ge 1).\endaligned$$

\par For a prime $p=24k+1=c^2+d^2=x^2+3y^2$ with $k,c,d,x,y\in\Bbb
Z$ and $c\e 1\mod 4$, in [HW] and [H], by using cyclotomic numbers
and Jacobi sums Hudson and Williams proved that
$$3^{\f{p-1}8}\e \cases \pm 1\mod p&\t{if $c\e \pm (-1)^{\f y4}\mod
3$,}
\\\pm \f dc\mod p&\t{if $d\e \pm (-1)^{\f y4}\mod 3$.}
\endcases$$
Let $p$ be a prime of the form $4k+1$, $q\in\Bbb Z$, $2\nmid q$ and
$p\nmid q$. Suppose that $p=c^2+d^2=x^2+qy^2$ with $c,d,x,y\in\Bbb
Z$ and $c\e 1\mod 4$. In [S5] and [S6] the author posed many
conjectures on $q^{[p/8]}$ and $\sls{b+\sqrt
{a^2+b^2}}2^{(p-1)/4}\mod p$ in terms of $c,d,x$ and $y$, where
$[\cdot]$ is the greatest integer function. For $m,n\in\Bbb Z$ let
$(m,n)$ be the greatest common divisor of $m$ and $n$. For $m\in\Bbb
Z$ with $m=2^{\alpha}m_0(2\nmid m_0)$ we say that $2^{\alpha} \
\Vert\ m$. In the paper, by developing the calculation technique of
quartic Jacobi symbols we partially solve many conjectures in [S5]
and [S6], and establish new reciprocity laws for quartic and octic
residues on condition that $(c,x+d)=1$ or $(d,x+c)=2^{\a}$. For the
history of classical reciprocity laws, see [Lem]. Suppose
$d=2^rd_0$, $y=2^ty_0$ and $d_0\e y_0\e 1\mod 4$. Assume $(c,x+d)=1$
or $(d_0,x+c)=1$. We then have the following typical results.

\par (1.1) If $p\e q\e 1\mod 8$, $q$ is a prime
and $q=a^2+b^2$ with $a,b\in\Bbb Z$, then
$$q^{\f{p-1}8}\e (-1)^{\f d4+\f{xy}4}\Ls dc^m\mod p
\iff \Ls{ac+bd}{ac-bd}^{\f{q-1}8}\e \Ls ba^m\mod q .$$
\par (1.2) If $q\e 7\mod 8$ is a prime, then
$$\aligned&q^{[p/8]}\e\cases (-1)^{\f y4}\sls dc^m\mod p&\t{if $p\e 1\mod
8$,}\\(-1)^{\f{x-1}2}\sls dc^m\f yx\mod p&\t{if $p\e 5\mod 8$.}
\endcases\\&\iff \Ls{c-di}{c+di}^{\f{q+1}8}\e i^m\mod q.\endaligned$$
\par (1.3) If $p\e 1\mod 8$,
$q=a^2+b^2$, $a,b\in\Bbb Z$, $2\mid a$ and $(a,b)=1$, then
$$\aligned &q^{\f{p-1}8}\e\cases (-1)^{\f d4+\f x4}\sls cd^m
\mod p&\t{if $4\mid a$ and $2\mid x$,}
\\(-1)^{\f d4+\f y4}\sls cd^m\mod p&\t{if $4\mid a$ and $2\nmid x$,}
\\(-1)^{\f{b-1}2+\f d4+\f {x+2}4}\sls cd^{m-1}
\mod p&\t{if $2\ \Vert\ a$ and $2\mid x$,}
\\(-1)^{\f{b-1}2+\f d4+\f y4+\f{x-1}2}\sls cd^{m-1}\mod p
&\t{if $2\ \Vert\ a$ and $2\nmid x$}
\endcases\\&\iff \Qs{(ac+bd)/x}{b+ai}=i^m.\endaligned$$
\par (1.4) For $q=b^2+16$ we have
$$\Ls{b+\sqrt{b^2+16}}2^{\f{p-1}4}\e
(-1)^{\f{b-1}2\cdot\f{p-1}4+[\f{b-3}4]y}\mod p.$$
\par (1.5) For $p\e 1\mod 8$ and $q=b^2+4$ we have
$$p\mid U_{\f{p-1}8}(b,-1)\iff\f{p-1}8+\f d4+\f y4\e
0\mod 2.$$
\par (1.6) For $p\e 1\mod 8$ and $q=4a^2+1$ we have
$$p\mid U_{\f{p-1}8}(4a,-1)\iff \f{p-1}8\e\cases \f d4+\f
a2y+\f{xy}4\mod 2&\t{if $2\mid a$,}\\\f d4+\f y4\mod 2&\t{if $2\nmid
a$.}\endcases$$

 \subheading{2. Basic lemmas}

 \pro{Lemma 2.1 ([S5, Proposition 2.1])} Let $a,b\in\Bbb Z$ with $2\nmid a$
and $2\mid b$. Then
$$\Qs i{a+bi}=i^{\f{a^2+b^2-1}4}
=(-1)^{\f{a^2-1}8}i^{(1-(-1)^{\f b2})/2}$$ and
$$\aligned\Qs{1+i}{a+bi}=\cases i^{((-1)^{\f{a-1}2}(a-b)-1)/4}&\t{if $4\mid b$,}
\\i^{\f{(-1)^{\f{a-1}2}(b-a)-1}4-1}&\t{if $2\ \Vert\ b$.}\endcases
\endaligned$$
\endpro
\pro{Lemma 2.2 ([S5, Proposition 2.2])} Let $a,b\in\Bbb Z$ with
$2\nmid a$ and $2\mid b$. Then
$$\Qs{-1}{a+bi}=(-1)^{\f b2}\qtq{and}\Qs
2{a+bi}=i^{(-1)^{\f{a-1}2}\f b2}.$$
\endpro
\pro{Lemma 2.3 ([S5, Proposition 2.3])} Let $a,b,c,d\in\Bbb Z$ with
$2\nmid ac$, $2\mid b$ and $2\mid d$. If $a+bi$ and $c+di$ are
relatively prime elements of $\Bbb Z[i]$, then we have the following
general law of quartic reciprocity:
$$\Qs{a+bi}{c+di}=(-1)^{\f b2\cdot \f{c-1}2+\f d2\cdot \f{a+b-1}2}
\Qs{c+di}{a+bi}.$$ In particular, if $4\mid b$, then
$$\Qs{a+bi}{c+di}=(-1)^{\f{a-1}2\cdot\f d2}\Qs{c+di}{a+bi}.$$
\endpro
\pro{Lemma 2.4 ([E], [S2, Lemma 2.1])} Let $m\in\Bbb N$ and
$a,b\in\Bbb Z$ with $2\nmid m$ and $(m,a^2+b^2)=1$. Then
$$\Qs{a+bi}m^2=\Ls{a^2+b^2}m.$$
\endpro

 \pro{Lemma 2.5 ([S4, Lemma 4.3])} Let $a,b\in\Bbb Z$ with
$2\nmid a$ and $2\mid b$. For any integer $x$ with $(x,a^2+b^2)=1$
we have
$$\Qs{x^2}{a+bi}=\Ls x{a^2+b^2}.$$
\endpro

\pro{Lemma 2.6} Let $a,b\in\Bbb Z$ with $2\mid b$ and $(a,b)=1$.
Then
$$\Qs b{a+bi}=\cases 1&\t{if $4\mid b$,}
\\(-1)^{\f{a-1}2}i&\t{if $2\ \Vert\ b$.}
\endcases$$
\endpro
Proof. By Lemmas 2.1 and 2.3 we have
$$\aligned\Qs b{a+bi}&=\Qs i{a+bi}\Qs{-bi}{a+bi}
=\Qs i{a+bi}\Qs a{a+bi}
\\&=\Qs i{a+bi}\cdot (-1)^{\f{a-1}2\cdot \f b2}\Qs{a+bi}a
=(-1)^{\f{a-1}2\cdot \f b2}\Qs i{a+bi}\Qs ia
\\&=(-1)^{\f{a-1}2\cdot \f b2}\cdot
(-1)^{\f{a^2-1}8}i^{\f{1-(-1)^{b/2}}2}\cdot (-1)^{\f{a^2-1}8}
\\&=\cases 1&\t{if $4\mid b$,}
\\(-1)^{\f{a-1}2}i&\t{if $2\ \Vert\ b$.}
\endcases\endaligned$$
Thus the lemma is proved.

\par For a given odd prime $p$ let $\Bbb Z_p$ denote the set of
 those rational numbers whose denominator is not divisible by $p$.
 Following [S2,S3] we define
 $$Q_r(p)=\Big\{k\bigm|\ k\in\Bbb Z_p,\
 \Qs{k+i}p=i^r\Big\}\qtq{for}r=0,1,2,3.$$
\pro{Lemma 2.7 ([S2, Theorem 2.3])} Let $p$ be an odd prime,
$r\in\{0,1,2,3\}$, $k\in\Bbb Z_p$ and $k^2+1\not\e 0\mod p$.
\par $(\t{\rm i})$ If $p\e 1\mod 4$ and $t^2\e -1\mod p$ with
 $t\in\Bbb Z_p$, then $k\in Q_r(p)$ if and only if
 $\sls{k+t}{k-t}^{(p-1)/4}\e t^r\mod p$.
\par $(\t{\rm ii})$ If $p\e 3\mod 4$, then $k\in Q_r(p)$ if and only if
 $\sls{k-i}{k+i}^{(p+1)/4}\e i^r\mod p$.\endpro
\pro{Lemma 2.8} Let $p$ be an odd prime, $k\in\Bbb Z_p$ and $n^2\e
k^2+1\mod p$ with $n\in\Bbb Z_p$ and $n(n+1)\not\e 0\mod p$. Then
$\sqs{k+i}p=\sls{n(n+1)}p$.\endpro Proof. For $k\e 0\mod p$ we have
$\sqs{k+i}p=\sqs ip=(-1)^{\f{p^2-1}8}=\sls{1\cdot 2}p$. So the
result is true. Now assume $k\not\e 0\mod p$. Then
$\sls{n-1}p\sls{n+1}p=\sls{n^2-1}p=\sls{k^2}p=1$ and so
$\sls{n-1}p=\sls{n+1}p$. By Lemma 2.4, $\sqs{k+i}p^2=\sls{k^2+1}p=1$
and so $\sqs{k+i}p=\pm 1$. By [S2, Theorem 2.4], $\sqs{k+i}p=1\iff
k\in Q_0(p)\iff \sls{n(n+1)}p=1$. Hence $\sqs{k+i}p=\sls{n(n+1)}p$.

\pro{Lemma 2.9} Suppose $c,d,m,x\in\Bbb Z$, $2\nmid m$, $x^2\e
c^2+d^2\mod m$ and $(m,x(x+d))=1$. Then
$$\Qs{c+di}m=\Ls{x(x+d)}m.$$
\endpro
Proof. Suppose that $p$ is a prime divisor of $m$. Then $p\nmid
x(x+d)$. If $p\nmid d$, then $\sls xd^2\e \sls cd^2+1\mod p$. Thus,
applying Lemma 2.8 we obtain
$$\Qs{c+di}p=\Qs{\f cd+i}p=\Ls{\f xd(1+\f xd)}p=\Ls{x(x+d)}p.$$
 When $p\mid
d$, we  have $p\nmid c$ and so
$$\Qs{c+di}p=\Qs cp=1=\Ls{x^2}p=\Ls{x(x+d)}p.$$
Hence,
$$\Qs{c+di}m=\prod_{p\mid m}\Qs{c+di}p=\prod_{p\mid
m}\Ls{x(x+d)}p=\Ls{x(x+d)}m,$$ where in the product $p$ runs over
all prime divisors of $m$. The proof is now complete.

\pro{Lemma 2.10} Suppose $c,d,x,y,q\in\Bbb Z$, $2\nmid c$, $2\mid
d$,
  $c^2+d^2=x^2+qy^2$, $y=2^ty_0$, $y_0 \e 1\mod 4$ and
$(y_0,x(x+d))=1$.
\par $(\t{\rm i})$ If $2\mid x$, then
$$\Ls{y}{c^2+(x+d)^2}=(-1)^{\f{y-1}4}
\Qs{y^{-1}}{c+di}.$$
\par $(\t{\rm ii})$ If $2\nmid x$, then
$$\Ls{y}{(c^2+(x+d)^2)/2}=(-1)^{\f{c^2-(x+d)^2}8t}
i^{(-1)^{\f{c-1}2}\f d2t}\Qs{y^{-1}}{c+di}.$$
\endpro
Proof. Since $c^2+(x+d)^2=2x(x+d)+qy^2$ we see that
$(c^2+(x+d)^2,y_0)=1$. For even $x$ we have $2\nmid qy$,
$c^2+(x+d)^2\e 1\mod 4$ and so
$$\Ls{y}{c^2+(x+d)^2}=\Ls{c^2+(x+d)^2}y
=\Ls {2x(x+d)}y =(-1)^{\f{y-1}4}\Ls{x(x+d)}y.$$ For odd $x$ we have
$c^2+(x+d)^2\e 2\mod 8$ and so
$$\align &\Ls{y}{(c^2+(x+d)^2)/2}
\\&=\Ls{2^ty_0}{(c^2+(x+d)^2)/2}
=(-1)^{\f{(c^2+(x+d)^2)/2-1}4t}\Ls{(c^2+(x+d)^2)/2}{y_0}
\\&=(-1)^{\f{c^2-(x+d)^2}8t}\Ls 2{y_0}\Ls {2x(x+d)+qy^2}{y_0}
=(-1)^{\f{c^2-(x+d)^2}8t}\Ls{x(x+d)}{y_0}.\endalign$$
 Since $x^2\e
c^2+d^2\mod{y_0}$, using Lemmas 2.9, 2.3 and 2.2  we see that
$$\align\Ls {x(x+d)}{y_0}
&=\Qs{c+di}{y_0}=\Qs{y_0}{c+di}=\Qs{y_0^{-1}}{c+di} \\&=\Qs
2{c+di}^t\Qs{y^{-1}}{c+di}=i^{(-1)^{\f{c-1}2}\f
d2t}\Qs{y^{-1}}{c+di}.\endalign$$ Now combining all the above we
obtain the result.

 \pro{Lemma
2.11} Suppose $c,d,x,y,q\in\Bbb Z$, $c\e 1\mod 4$, $2\mid d$,
  $c^2+d^2=x^2+qy^2$, $y=2^ty_0$, $y_0 \e 1\mod 4$,
$(y_0,x(x+c))=1$ and $2^n\ \Vert\ ((x+c)^2+d^2)$. Then
$$\Ls{y}{((x+c)^2+d^2)/2^n}=(-1)^{\f{\f{(x+c)^2+d^2}{2^n}-1}4t
+\f{y_0-1}4n}i^{\f d2t}\Qs{y^{-1}}{c+di}.$$
\endpro
Proof. It is well known that the sum of two squares has no prime
divisors of the form $4k+3$. Hence $((x+c)^2+d^2)/2^n\e 1\mod 4$ and
so
$$\align &\Ls{y}{((x+c)^2+d^2)/2^n}\\&=
\Ls{2^ty_0}{((x+c)^2+d^2)/2^n}=(-1)^{\f{\f{(x+c)^2+d^2}{2^n}-1}4t}
\Ls {((x+c)^2+d^2)/2^n}{y_0}
\\&=(-1)^{\f{\f{(x+c)^2+d^2}{2^n}-1}4t}\cdot (-1)^{\f{y_0-1}4n}
\Ls {(x+c)^2+d^2}{y_0} .\endalign$$ Since $x^2\e c^2+d^2\mod{y_0}$,
using Lemma 2.9 we see that
$$\align \Ls{(x+c)^2+d^2}{y_0}&
=\Ls {2x(x+c)}{y_0}=\Ls 2{y_0}\Qs{-d+ci}{y_0}
\\&=\Qs{-i}{y_0}\Qs{-d+ci}{y_0}=\Qs{c+di}{y_0}=\Qs{y_0}{c+di}
\\&=\Qs{y_0^{-1}}{c+di}=\Qs{2^ty^{-1}}{c+di}
=i^{\f d2t}\Qs{y^{-1}}{c+di}.\endalign$$ Now combining all the above
we obtain the result.

\pro{Lemma 2.12} Let $p$ be  a prime of the form $4k+1$ and
$p=c^2+d^2$ with $c,d\in\Bbb Z$. Suppose $q\in\Bbb Z$, $p\nmid q$
and $p=x^2+qy^2$ with $x,y\in\Bbb Z$. Then $(x+d,c^2)=(x+d,qy^2)$
and
$$(qy^2,c^2+(x+d)^2)=(x+d,c^2)\Big(2,x+d+\f{c^2}{(x+d,c^2)}\Big).$$
\endpro
Proof.  Since $(x,y)^2\mid p$ we see that $(x,y)=1$. If $p\mid x$,
then $p\mid qy^2$ and so $p\mid y$. This contradicts the fact
$(x,y)=1$. Hence $p\nmid x$. Since
$(x,c^2+(x+d)^2)=(x,c^2+d^2)=(x,p)=1$ and
$qy^2=d^2-x^2+c^2=c^2+(x+d)^2-2x(x+d),$
 we see that $(x+d,c^2)=(x+d,x^2-d^2+qy^2)=(x+d,qy^2)$ and
$$\align &(qy^2, c^2+(x+d)^2)\\&=(2x(x+d),c^2+(x+d)^2)
=(2(x+d),c^2+(x+d)^2)
\\&=(x+d,c^2)\Big(2\f{x+d}{(x+d,c^2)},\f{c^2}{(x+d,c^2)}+(x+d)
\f{x+d}{(x+d,c^2)}\Big)
\\&=(x+d,c^2)\Big(2,\f{c^2}{(x+d,c^2)}+(x+d)
\f{x+d}{(x+d,c^2)}\Big)
\\&=(x+d,c^2)\Big(2,\f{c^2}{(x+d,c^2)}+(x+d)\Big).
\endalign$$
 Thus  the lemma is proved.

 \pro{Lemma 2.13} Let $p$ be  a prime of the
form $4k+1$ and $p=c^2+d^2$ with $c,d\in\Bbb Z$ and $2\nmid c$.
Suppose $q\in\Bbb Z$, $p\nmid q$, $p=x^2+qy^2$, $x,y\in\Bbb Z$ and
$\sqs{x/y}{c+di}=(-1)^{[\f p8]+n}i^k$. Then
$$\aligned q^{[p/8]}\e \cases (-1)^{n}\sls dc^k\mod p&\t{if $p\e 1\mod 8$,}
\\ (-1)^{n}\sls dc^k\f yx\mod p&\t{if $p\e 5\mod
8.$}\endcases\endaligned$$
\endpro
Proof. It is clear that $(c,d)=1$, $p\nmid y$ and so
$$\Ls xy^{\f{p-1}4}\e \Qs{x/y}{c+di}=(-1)^{[\f p8]+n}i^k\e
(-1)^{[\f p8]+n}\Ls dc^k \mod {c+di}.$$ Thus
$$\Ls xy^{\f{p-1}4}\e (-1)^{[\f p8]+n}\Ls dc^k
\mod p.$$ Hence, for $p\e 1\mod 8$ we have
$$q^{\f{p-1}8}=(-1)^{\f{p-1}8}(-q)^{\f{p-1}8}
\e (-1)^{\f{p-1}8}\Ls xy^{\f{p-1}4}\e (-1)^n\Ls dc^k \mod p,$$ for
$p\e 5\mod 8$ we have
$$q^{\f{p-5}8}=(-1)^{\f{p-5}8}(-q)^{\f{p-5}8}
\e (-1)^{\f{p-5}8}\Ls xy^{\f{p-5}4}\e (-1)^n\Ls dc^k \f yx\mod p.$$
This proves the lemma.

\subheading{3. Determination of $q^{[p/8]}\mod p$ using
$\sqs{c+(x+d)i}q$ or $\sqs{d-(x+c)i}q$}
 \pro{Theorem 3.1} Let $p$ be a prime of the form
$4k+1$, $q\in\Bbb Z$, $2\nmid q$ and $p\nmid q$. Suppose that
$p=c^2+d^2=x^2+qy^2$ with $c,d,x,y\in\Bbb Z$, $c\e  1\mod 4$,
$d=2^rd_0$, $d_0\e 1\mod 4$, $(c,x+d)=1$, $2\mid x$, $y\e 1\mod 4$
and $\sqs{c/(x+d)+i}q=i^k$.
\par $(\t{\rm i})$ If $p\e 1\mod 8$, then
$$q^{\f{p-1}8}\e
\cases (-1)^{\f{q-1}8+\f d4+\f x4}\sls dc^k\mod p&\t{if $q\e 1\mod
8$,}\\(-1)^{\f{q-5}8+\f d4+\f {x-2}4}\sls dc^{k+1}\mod p&\t{if $q\e
5\mod 8$.}
\endcases$$

\par $(\t{\rm ii})$ If $p\e 5\mod 8$, then
$$q^{\f{p-5}8}\e
\cases (-1)^{\f{q-1}8+\f {x-2}4}\sls dc^{k+1}\f yx\mod p&\t{if $q\e
1\mod 8$,}\\(-1)^{\f{q-5}8+\f x4}\sls dc^k\f yx\mod p&\t{if $q\e
5\mod 8$.}
\endcases$$

  \endpro

Proof. Suppose $2^m\ \Vert\ (x+d)$ and $x=2^sx_0(2\nmid x_0)$. Since
$2\mid x$ we have $q\e 1\mod 4$. As $(c,x+d)=1$, by Lemma 2.12 we
have $(qy,x+d)=1$ and $(qy^2,c^2+(x+d)^2)=1$. Note that $(x,y)^2\mid
p$. We also have $(x,y)=1$.   Using Lemmas 2.1-2.5, 2.10 and the
fact that $\sqs an=1$ for $a,n\in\Bbb Z$ with $2\nmid n$ and
$(a,n)=1$ we see that
$$\align i^k&=\Qs{c+(x+d)i}q
=\Qs q{c+(x+d)i} =\Qs {qy^2}{c+(x+d)i}\Qs{y^2}{c+(x+d)i}
\\&=
\Qs{c^2+d^2-x^2}{c+(x+d)i} \Ls
y{c^2+(x+d)^2}
\\&=
\Qs{-2x(x+d)+c^2+(x+d)^2}{c+(x+d)i} \Ls y{c^2+(x+d)^2}
\\&=
\Qs{-2x(x+d)}{c+(x+d)i}(-1)^{\f{y-1}4} \Qs{y^{-1}}{c+di}
\endalign$$ and
$$\align &\Qs{-2x(x+d)}{c+(x+d)i} \\&=\Qs 2{c+(x+d)i}^{m+s+1}
\Qs{-x_0(x+d)/2^m}{c+(x+d)i}
\\&=i^{\f{x+d}2(m+s+1)}(-1)^{\f{x_0(x+d)/2^m+1}2\cdot\f{x+d}2}\Qs
{c+(x+d)i}{x_0(x+d)/2^m}
\\&=(-1)^{\f{x_0(x+d)/2^m+1}2\cdot\f{x+d}2}
i^{\f{x+d}2(m+s+1)}\Qs{c+di}{x_0}\Qs c{(x+d)/2^m}
\\&=(-1)^{\f{x_0(x+d)/2^m+1}2\cdot\f{x+d}2}
i^{\f{x+d}2(m+s+1)}(-1)^{\f{x_0-1}2\cdot\f d2}\Qs{x_0}{c+di}
\\&=(-1)^{\f{x_0(x+d)/2^m+1}2\cdot\f{x+d}2}
i^{\f{x+d}2(m+s+1)}(-1)^{\f{x_0-1}2\cdot\f d2}\Qs
2{c+di}^{-s}\Qs{x}{c+di}
\\&=(-1)^{\f{x_0(x+d)/2^m+1}2\cdot\f{x+d}2}
i^{\f{x+d}2(m+s+1)}(-1)^{\f{x_0-1}2\cdot\f d2} i^{-\f
d2s}\Qs{x}{c+di}.\endalign$$
 Therefore,
$$\aligned i^k&=(-1)^{\f{x_0(x+d)/2^m+1}2\cdot\f{x+d}2+
\f{x_0-1}2\cdot \f d2+\f{y-1}4}i^{\f{x+d}2(m+s+1)-\f
d2s}\Qs{x/y}{c+di}.
\endaligned$$
Observe that
$$\aligned (-1)^{\f{y-1}4}
=(-1)^{\f{q(y^2-1)}8}=(-1)^{\f{p-q-x^2}8}=\cases
(-1)^{\f{p-q-4}8}&\t{if $2\ \Vert\ x$,}
\\(-1)^{\f{p-q}8}&\t{if $4\mid x$}\endcases\endaligned$$
 and
$$\aligned &i^{\f
d2s-\f{x+d}2(m+s+1)}=i^{-\f{x+d}2(m+1)-\f x2s}
=(-1)^{\f{(m+1)(x+d)}4}i^{-\f x2s}
 \\&=\cases (-1)^{\f{(m+1)(x+d)}4}\cdot (-1)^{\f{x+2}4}i&\t{if $2\
 \Vert\ x$,}
 \\(-1)^{\f{(m+1)(x+d)}4}&\t{if $4\mid x$.}\endcases\endaligned$$
 We then obtain
$$\aligned &\Qs{x/y}{c+di}\\&=\cases
(-1)^{\f{x_0(x+d)/2^m+1}2\cdot\f{x+d}2+ \f{x-2}4\cdot \f
d2+\f{p-q-4}8}\cdot(-1)^{\f{(m+1)(x+d)}4+\f{x+2}4}i^{k+1}&\t{if $2\
\Vert\ x$,}
\\(-1)^{\f{x_0(x+d)/2^m+1}2\cdot\f{x+d}2+ \f{x_0-1}2\cdot \f
d2+\f{p-q}8}\cdot(-1)^{\f{(m+1)(x+d)}4}i^k&\t{if $4\mid
x$.}\endcases\endaligned\tag 3.1$$

\par If $p\e q\e 1\mod 8$, then $4\mid d$,
$4\mid x$ and
$4\mid (x+d)$. By (3.1) we have
$$\Qs{x/y}{c+di}=(-1)^{\f{p-q}8+\f{x+d}4}i^k=(-1)^{\f{p-1}8+\f{q-1}8+\f
d4+\f x4}i^k.$$

\par If $p\e 1\mod 8$ and $q\e 5\mod 8$, then
$2\ \Vert\ x$, $4\mid d$, $2\ \Vert\ (x+d)$, $m=1$ and
$(-1)^{\f{x_0(x+d)/2+1}2}=(-1)^{\f d4x_0+\f{x_0^2+1}2}=(-1)^{\f
d4+1}$. Thus, by (3.1) we get
$$\align \Qs{x/y}{c+di}
&=(-1)^{\f{x_0(x+d)/2+1}2+\f{p-q-4}8+\f{x+d}2+\f{x+2}4}i^{k+1}
\\&=(-1)^{\f d4+1+\f{p-1}8-\f{q-5}8+\f{x+2}4}i^{k+1}
=(-1)^{\f{p-1}8+\f{q-5}8+\f d4+\f{x-2}4}i^{k+1}.
\endalign$$
\par If $p\e 5\mod 8$ and $q\e 1\mod 8$, then
$2\ \Vert\ d$, $2\ \Vert\ x$, $4\mid (x+d)$ and $m\ge 2$. From (3.1)
we deduce that
$$\align \Qs{x/y}{c+di}
&=(-1)^{\f{x-2}4+\f{p-q-4}8+\f{x+d}4+\f{x+2}4}i^{k+1}
\\&=(-1)^{\f{x-2}4+\f {p-5}8-\f{q-1}8+\f{4x_0+2d_0+2}4}i^{k+1}
=(-1)^{\f{p-5}8+\f{q-1}8+\f{x-2}4}i^{k+1}.
\endalign$$
\par If $p\e q\e 5\mod 8$, then $2\ \Vert\ d$,
$4\mid x$, $2\ \Vert\ (x+d)$, $m=1$ and
$(-1)^{\f{x_0(x+d)/2+1}2}=(-1)^{\f x4x_0+\f{d_0x_0+1}2}=(-1)^{\f
x4+\f{x_0+1}2}$. Thus, by (3.1) we have
$$\align \Qs{x/y}{c+di}
&=(-1)^{\f{x_0(x+d)/2+1}2+\f{x_0-1}2+\f{p-q}8 +\f{x+d}2}i^k
\\&=(-1)^{\f x4+\f{x_0+1}2+\f{x_0-1}2+\f{p-q}8+1}i^k
=(-1)^{\f{p-5}8+\f{q-5}8+\f x4}i^k.
\endalign$$
\par Now combining the above with Lemma 2.13 we deduce the result.
\pro{Lemma 3.1} Let $p$ be a prime of the form $4k+1$, $q\in\Bbb Z$,
$2\nmid q$ and $p\nmid q$. Suppose that $p=c^2+d^2=x^2+qy^2$ with
$c,d,x,y\in\Bbb Z$, $c\e  1\mod 4$, $d=2^rd_0$, $y=2^ty_0$, $d_0\e
y_0\e 1\mod 4$, $(c,x+d)=1$ and $2\nmid x$. Assume that
$\sqs{c/(x+d)+i}q=i^k$. Then
$$\aligned&\Qs{x/y}{c+di}
=\cases (-1)^{[\f p8]+\f{q-\sls{-1}q}4\cdot\f{x-1}2 + \f{q-1}2\cdot
\f d4}i^{\f{1-\sls{-1}qq}4+\f {(-1)^{\f{x-1}2}x-c}4+k}&\t{if $8\mid
(p-1)$,}
\\(-1)^{[\f p8]+\f{q+1}2+\f{x+1}2(1+\f{q-\sls{-1}q}4)} i^{\f{1-\sls{-1}qq}4
+\f{(-1)^{\f{x-1}2}x-c}4+k-1}
 &\t{if $8\mid (p-5)$.}
\endcases\endaligned$$
\endpro
Proof.  As $(c,x+d)=1$, by Lemma 2.11 we have $(qy_0,x+d)=1$ and
$(qy_0^2,c^2+(x+d)^2)=1$. Note that $(x,y)^2\mid p$. We also have
$(x,y)=1$. It is easily seen that
$$c+(x+d)i=i^{\f{1\mp 1}2}(1+i)\Big(\f{x+d\pm
c}2+\f{\pm(x+d)-c}2i\Big)$$ and so $$\Ls{x+d\pm
c}2^2+\Ls{\pm(x+d)-c}2^2=\f {c^2+(x+d)^2}2.$$ Set
$\ep=(-1)^{\f{p-1}4+\f{x-1}2}$. Since $4\mid d\iff 8\mid (p-1)$ we
see that $x+d\e \ep\mod 4$ and $4\mid (\ep(x+d)-c)$. Using Lemmas
2.1-2.5, 2.10 and the above we see that
$$\align i^k&
=\Qs{c+(x+d)i}q=\Qs iq^{\f{1-\ep}2}\Qs{1+i}q\Qs{\f{x+d+\ep
c}2+\f{\ep(x+d)-c}2i}q
\\&=(-1)^{\f{q-\sls{-1}q}4\cdot\f{1-\ep}2}i^{\f{\sls{-1}qq-1}4}
(-1)^{\f{q-1}2\cdot\f{\ep(x+d)-c}4}\Qs q{\f{x+d+\ep
c}2+\f{\ep(x+d)-c}2i}
\endalign$$
and
$$\align &\Qs q{\f{x+d+\ep
c}2+\f{\ep(x+d)-c}2i}
\\&=\Qs {qy^2}{\f{x+d+\ep
c}2+\f{\ep(x+d)-c}2i}\Qs {y^2}{\f{x+d+\ep c}2+\f{\ep(x+d)-c}2i}
\\&=\Qs {c^2+(x+d)^2-2x(x+d)}{\f{x+d+\ep
c}2+\f{\ep(x+d)-c}2i}\Ls y{\sls{x+d+\ep c}2^2+\sls{\ep(x+d)-c}2^2}
\\&=\Qs {2}{\f{x+d+\ep
c}2+\f{\ep(x+d)-c}2i}\Qs {-x(x+d)}{\f{x+d+\ep
c}2+\f{\ep(x+d)-c}2i}\Ls y{(c^2+(x+d)^2)/2}
\\&=i^{(-1)^{\f{(x+d+\ep c)/2-1}2}\f{\ep(x+d)-c}4}
(-1)^{\f{x(x+d)+1}2\cdot\f{\ep(x+d)-c}4}\Qs {\f{x+d+\ep
c}2+\f{\ep(x+d)-c}2i}{x(x+d)}
\\&\qq\times (-1)^{\f{c^2-(x+d)^2}8t}
i^{\f d2t}\Qs{y^{-1}}{c+di}.
\endalign$$
Obviously $i^{(-1)^ab}=(-1)^{ab}i^b$, $(-1)^{\f{x+d+\ep
c-2}4}=(-1)^{\f{\ep (x+d)+c-2\ep}4} =(-1)^{\f{\ep (x+d)-c}4+\f{1-\ep
}2}$ and so
$$\align i^{(-1)^{\f{x+d+\ep c-2}4}\f{\ep(x+d)-c}4}
&=(-1)^{(\f{\ep (x+d)-c}4+\f{1-\ep}2)\f{\ep (x+d)-c}4}i^{\f{\ep
(x+d)-c}4}\\& =(-1)^{\f{1+\ep}2\cdot \f{\ep (x+d)-c}4}i^{\f{\ep
(x+d)-c}4}.\endalign$$ Also,
$(-1)^{\f{x(x+d)+1}2}=(-1)^{\f{x^2-1}2+\f d2x+1}=(-1)^{\f d2+1}$ and
$(-1)^{\f{c^2-(x+d)^2}8}=(-1)^{\f{\ep (x+d)-c}4}$. Thus,
$$\aligned &i^{(-1)^{\f{(x+d+\ep c)/2-1}2}\f{\ep(x+d)-c}4}
(-1)^{\f{x(x+d)+1}2\cdot\f{\ep(x+d)-c}4}\cdot
(-1)^{\f{c^2-(x+d)^2}8t} i^{\f d2t}
\\&=(-1)^{\f{1+\ep}2\cdot \f{\ep(x+d)-c}4}i^{\f{\ep(x+d)-c}4}
(-1)^{(\f d2+1)\f{\ep(x+d)-c}4}\cdot (-1)^{\f{\ep(x+d)-c}4t}i^{\f
d2t}
\\&=(-1)^{(\f{1-\ep}2+\f
d2+t)\f{\ep(x+d)-c}4}i^{\f{\ep(x+d)-c}4+\f d2t}.
\endaligned$$
It is easily seen that
$$\align&\Qs {\f{x+d+\ep
c}2+\f{\ep(x+d)-c}2i}{x(x+d)}
\\&=\Qs {x+d+\ep
c+(\ep(x+d)-c)i}{x}\Qs {x+d+\ep c+(\ep(x+d)-c)i}{x+d}
\\&=\Qs{d+\ep c+(\ep d-c)i}x\Qs{\ep c-ci}{x+d}=\Qs{(\ep-i)(c+di)}x
\Qs{\ep -i}{x+d}
\\&=\Qs{\ep -i}{x(x+d)}\Qs{c+di}x=
\Qs{i^{\f{5+\ep}2}(1+i)}{x(x+d)}\Qs{c+di}x
\\&=\Qs i{x(x+d)}^{\f{5+\ep}2}\Qs{1+i}{x(x+d)}(-1)^{\f{x-1}2\cdot\f
d2}\Qs x{c+di}
\\&=(-1)^{\f{(x(d+x))^2-1}8\cdot\f{5+\ep}2}
i^{\f{(-1)^{d/2}x(x+d)-1}4} (-1)^{\f{x-1}2\cdot\f d2}\Qs
x{c+di}.\endalign$$  Now combining all the above we deduce that
$$\aligned i^k
&=(-1)^{\f{q-\sls{-1}q}4\cdot\f{1-\ep}2
+\f{q-1}2\cdot\f{\varepsilon(x+d)-c}4} i^{\f{\sls{-1}qq-1}4}
\\&\q\times (-1)^{(\f{1-\varepsilon}2+\f
d2+t)\f{\varepsilon(x+d)-c}4}i^{\f{\varepsilon(x+d)-c}4+\f d2t}
\\&\q\times(-1)^{\f{(x(d+x))^2-1}8\cdot\f{5+\varepsilon}2
+\f{x-1}2\cdot \f d2} i^{\f{(-1)^{d/2}x(x+d)-1}4} \Qs{x/y}{c+di}.
\endaligned\tag 3.2$$
It is clear that
$$\aligned (-1)^{\f{\ep(x+d)-c}4}&=(-1)^{\f{\ep(x+d)-c}4\cdot \f{\ep(x+d)+c}2}
=(-1)^{\f{(x+d)^2-c^2}8}=(-1)^{\f{2d^2+2dx-qy^2}8}
\\&=\cases (-1)^{\f d4}&\t{if $8\mid (p-1)$ and so $d\e y\e 0\mod
4$,}
\\(-1)^{\f{2+x-q}2}=(-1)^{\f{q-1}2+\f{1-\ep}2}&\t{if $8\mid (p-5)$
and so $d\e y\e 2\mod 4$}
\endcases\endaligned$$
and
$$\aligned (-1)^{\f{x^2(x+d)^2-1}8}&=(-1)^{\f {(-1)^{d/2}x(x+d)-1}4}
=(-1)^{\f{(-1)^{d/2}(dx+1)-1}4}
\\&=\cases (-1)^{\f d4}&\t{if $8\mid (p-1)$ and so $4\mid d$,}
\\(-1)^{\f{d_0x+1}2}=\ep&\t{if $8\mid (p-5)$
and so $d\e y\e 2\mod 4$.}
\endcases\endaligned$$
Since $x+d\e \ep\mod 4$ and $(-1)^{\f d2}=(-1)^{\f{p-1}4}$ we also
have
$$\aligned &i^{\f{\ep(x+d)-c}4+\f{(-1)^{d/2}x(x+d)-1}4+\f d2t}
\\&=(-1)^{\f{(-1)^{d/2}x(x+d)-1}4}i^{\f{\ep(x+d)-c}4-\f{(-1)^{(p-1)/4}x(x+d)-1}4+\f d2t}
\\&=(-1)^{\f{(-1)^{d/2}x(x+d)-1}4}i^{\ep(x+d)\f{1-(-1)^{(x-1)/2}x}4-\f{c-1}4+\f
d2t}
\\&=(-1)^{\f{(-1)^{d/2}x(x+d)-1}4+\f{c-1}4}i^{\f{1-(-1)^{(x-1)/2}x}4+\f{c-1}4+\f
d2t}
\\&=\cases (-1)^{\f d4+\f{c-1}4+\f d4t}i^{\f{c-(-1)^{(x-1)/2}x}4}
&\t{if $8\mid (p-1)$,}
\\\ep(-1)^{\f{c-1}4}i^{\f{c-(-1)^{(x-1)/2}x}4+d_0}=(-1)^{\f{1-\ep}2+\f{c-1}4}
i^{\f{c-(-1)^{(x-1)/2}x}4+1}&\t{if $8\mid (p-5)$.}
\endcases
\endaligned$$
Note that
$(-1)^{\f{c-1}4}=(-1)^{\f{c^2-1}8}=(-1)^{\f{p-1-d^2}8}=(-1)^{[\f
p8]}$. From the above and (3.2) we deduce the result.

\pro{Theorem 3.2} Let $p$ be a prime of the form $4k+1$, $q\in\Bbb
Z$, $2\nmid q$ and $p\nmid q$. Suppose that $p=c^2+d^2=x^2+qy^2$
with $c,d,x,y\in\Bbb Z$, $c\e  1\mod 4$, $d=2^rd_0$, $y=2^ty_0$,
$d_0\e y_0\e 1\mod 4$, $(c,x+d)=1$, $2\nmid x$ and
$\sqs{c/(x+d)+i}q=i^k$.
\par $(\t{\rm i})$ If $p\e 1\mod 8$, then
$$q^{\f{p-1}8}\e
\cases (-1)^{\f {q-1}8+\f d4+\f y4} \sls dc^k\mod p&\t{if $q\e 1\mod
8$,}
\\(-1)^{\f {q+5}8+\f{x-1}2+\f
y4}\sls dc^{k-1}\mod p
 &\t{if $q\e  3\mod 8$,}
\\(-1)^{\f {q-5}8+ \f d4+\f{x-1}2+\f
y4}\sls dc^{k-1}\mod p
 &\t{if $q\e  5\mod 8$,}
\\ (-1)^{\f {q+1}8+\f y4} \sls dc^k\mod p&\t{if $q\e 7\mod 8$.}
 \endcases$$
\par $(\t{\rm ii})$ If $p\e 5\mod 8$, then
$$q^{\f{p-5}8}\e
\cases (-1)^{\f{q-1}8+\f{x-1}2}\sls dc^{k-1}\f yx\mod p&\t{if $q\e
1\mod 8$,}\\(-1)^{\f{q+5}8}\sls dc^{k-1}\f yx\mod p&\t{if $q\e 3\mod
8$,}
\\(-1)^{\f{q+3}8}\sls dc^k\f yx\mod p&\t{if $q\e 5\mod
8$,} \\(-1)^{\f{q+1}8+\f{x-1}2}\sls dc^k\f yx\mod p&\t{if $q\e 7\mod
8$.}
\endcases$$

  \endpro
Proof. Suppose $2^m\ \Vert\ (c-(-1)^{(x-1)/2}x)$. Then $m\ge 2$ and
$2^{m+1}\ \Vert\ (c-x)(c+x)$. As $d^2-qy^2=-(c-x)(c+x)$ we see that
$2^{m+1}\ \Vert\ (d^2-qy^2)$. We first assume $p\e 1\mod 8$. Since
$4\mid d$ and $4\mid y$ we have $m\ge 3$ and $2^{m-3}\ \Vert\ (\sls
d4^2-q\sls y4^2)$. Thus,
$$ i^{\f{(-1)^{(x-1)/2}x-c}4}=(-1)^{\f{(-1)^{(x-1)/2}
x-c}8}=(-1)^{2^{m-3}}=(-1)^{\sls d4^2-q\sls y4^2} =(-1)^{\f d4+\f
y4}.$$ From the above and Lemma 3.1 we deduce that
$$\aligned \Qs{x/y}{c+di}=
\cases (-1)^{\f{p-1}8+\f {q-\sls{-1}q}8+\f{q+1}2\cdot \f d4+\f
y4}i^k &\t{if $q\e \pm 1\mod 8$,}
\\(-1)^{\f{p-1}8+\f {q-5\sls{-1}q}8+\f{q+1}2\cdot \f d4+\f{x-1}2+\f
y4}i^{k-1} &\t{if $q\e \pm 3\mod 8$.}\endcases
\endaligned$$
Now applying Lemma 2.12 we deduce (i).
\par Suppose $p\e 5\mod 8$. As
  $qy^2-d^2=4(qy_0^2-d_0^2)$ we get
  $2^{m-1}\ \Vert\ (qy_0^2-d_0^2)$. Clearly $qy_0^2-d_0^2\e q-1\mod
  8$. Thus
  $$q\e \cases 3\mod 4&\t{if $m=2$,}\\5\mod 8&\t{if $m=3$,}
  \\1\mod 8&\t{if $m>3$.}\endcases$$
 For $q\e 1\mod 4$ we have $8\mid (c-(-1)^{\f{x-1}2}x)$ and
 $$i^{\f{(-1)^{(x-1)/2}x-c}4}=(-1)^{\f{(-1)^{(x-1)/2}x-c}8}
 =(-1)^{2^{m-3}}=(-1)^{\f{q-1}4}.$$ Thus, using Lemma 3.1 we deduce that
$$\aligned \Qs{x/y}{c+di} =\cases
(-1)^{\f{p-5}8+\f{q-1}8+\f{x-1}2}i^{k-1} &\t{if $q\e 1\mod 8$,}
\\(-1)^{\f{p-5}8+\f{q+3}8}i^k &\t{if $q\e
5\mod 8$.}\endcases\endaligned$$ Now applying Lemma 2.12 we deduce
the result in the case $p\e 5\mod 8$ and $q\e 1\mod 4$.
 For  $q\e 3\mod 4$ we have $2^2\ \Vert\ (c-(-1)^{(x-1)/2}
x)$ and so
$$q\e qy_0^2=2c\cdot \f{c-(-1)^{\f{x-1}2} x}4-4\ls{c-(-1)^{\f{x-1}2} x}4^2+d_0^2
\e 2\cdot \f{c-(-1)^{\f{x-1}2}x}4-4+1\mod 8.$$ Thus,
$\f{(-1)^{(x-1)/2}x-c}4\e -\f{q+3}2\mod 4$ and so
$i^{\f{(-1)^{(x-1)/2}x-c}4}=i^{-\f{q+3}2}$. Using Lemma 3.1 we see
that
$$\aligned
\Qs{x/y}{c+di}=\cases (-1)^{\f{p-5}8+\f{q+5}8}i^{k-1} &\t{if $q\e
3\mod 8$,}
\\(-1)^{\f{p-5}8+\f{q+1}8+\f{x-1}2}i^k &\t{if $q\e
7\mod 8$.}\endcases\endaligned$$
 Now applying Lemma 2.12 we obtain
the result in the case $p\e 5\mod 8$ and $q\e 3\mod 4$.
\par Summarizing all the above we prove the theorem.

\par\q
\newline{\bf Remark 3.1}
We note that the $k$ in Theorems 3.1-3.2 depends only on $\f
c{x+d}\mod q$.

\pro{Corollary 3.1} Let $p\e 1,49\mod {60}$ be a prime and so
$p=c^2+d^2=x^2+15y^2$ with $c,d,x,y\in\Bbb Z$. Suppose $c\e 1\mod
4$, $d=2^rd_0$, $y=2^ty_0$, $d_0\e y_0\e 1\mod 4$ and $(c,x+d)=1$.
\par $(\t{\rm i})$ If $p\e 1\mod 8$, then
$$15^{\f{p-1}8}\e \cases (-1)^{\f y4}\mod p
&\t{if $\f c{x+d}\e 0,\pm 1\mod{15}$,}
\\ -(-1)^{\f y4}\mod p
&\t{if $\f c{x+d}\e \pm 4\mod{15}$,}
\\\pm(-1)^{\f y4}\f cd\mod p&\t{if $\f c{x+d}\e \pm 5,\pm 6\mod{15}$.}
\endcases$$
\par $(\t{\rm ii})$ If $p\e 5\mod 8$, then

$$15^{\f{p-5}8}\e \cases (-1)^{\f{x-1}2}\f yx\mod p
&\t{if $\f c{x+d}\e 0,\pm 1\mod{15}$,}
\\ -(-1)^{\f{x-1}2}\f yx\mod p
&\t{if $\f c{x+d}\e \pm 4\mod{15}$,}
\\\mp(-1)^{\f {x-1}2}\f {dy}{cx}
\mod p&\t{if $\f c{x+d}\e \pm 5,\pm 6\mod{15}$.}
\endcases$$
\endpro
Proof. Clearly $x$ is odd . Thus, putting $q=15$ in Theorem 3.2 and
noting that (see [S2, Example 2.1])
$$\Qs{n+i}{15}=\Qs{n+i}3\Qs{n+i}5=\cases 1&\t{if $n\e 0,\pm 1\mod{15}$,}
\\-1 &\t{if $n\e \pm 4\mod{15}$,}
\\\mp i &\t{if $n\e \pm 5,\pm 6\mod{15}$}
\endcases$$ we deduce the result.
\par\q
\par For example, since $61=5^2+(-6)^2=(-1)^2+15\cdot 2^2$, $(5,-1-6)=1$
 and
$\f 5{-1-6}\e -5\mod{15}$ we have
$$15^{\f{61-5}8}=15^7\e (-1)^{\f{-1-1}2}\f{-6\cdot 2}{5\cdot (-1)}
\e -\f{12}5\e 22\mod{61}.$$

\pro{Theorem 3.3} Let $p$ be a prime of the form $4k+1$, $q\in\Bbb
Z$, $2\nmid q$ and $p\nmid q$. Suppose that $p=c^2+d^2=x^2+qy^2$
with $c,d,x,y\in\Bbb Z$, $c\e 1\mod 4$, $d=2^rd_0$, $y=2^ty_0$,
$d_0\e y_0\e 1\mod 4$, $(d_0,x+c)=1$, $2^m\ \Vert\ (x+c)$ and
$\sqs{d/(x+c)-i}q=i^k$.
\par $(\t{\rm i})$ If $2\mid x$, then
$$q^{[p/8]}\e
\cases (-1)^{\f d4+[\f x4]}\sls dc^k\mod p&\t{if $p\e 1\mod 8$,}
\\(-1)^{[\f{x+2}4]}\sls dc^{k-1}\f yx\mod p&\t{if $p\e 5\mod 8$.}
\endcases$$

$(\t{\rm ii})$ If $x\e 1\mod 4$ and $p\e 1\mod 8$, then
$$q^{\f{p-1}8}\e (-1)^{\f{q^2-1}8+\f{q+1}2\cdot \f d4+\f y4}\Ls
dc^k\mod p.$$

$(\t{\rm iii})$ If $x\e 3\mod 4$, $p\e 1\mod 8$ and $m<r$, then
$$q^{\f{p-1}8}\e (-1)^{\f{q+1}2\cdot \f d4+\f y4}\Ls
dc^k\mod p.$$

  \endpro

Proof. Suppose  $x=2^sx_0(2\nmid x_0)$. As $(x+c,d_0)=1$, by Lemma
2.12 we have $(qy_0,x+c)=1$ and $(qy_0^2,(x+c)^2+d^2)=1$. Note that
$(x,y)^2\mid p$. We also have $(x,y)=1$. In all three cases we have
$m<r$.  Using Lemmas 2.1-2.5 and the fact that $\sqs aq=1$ for
$a\in\Bbb Z$ with $(a,q)=1$ we see that
$$\align \Qs{d-(x+c)i}q&=\Qs {-i}q\Qs{x+c+di}q
=(-1)^{\f{q^2-1}8}\Qs{\f{x+c}{2^m}+\f d{2^m}i}q
\\&=(-1)^{\f{q^2-1}8+\f{q-1}2\cdot\f d{2^{m+1}}}\Qs q{\f{x+c}{2^m}+\f
d{2^m}i}
\\&=(-1)^{\f{q^2-1}8+\f{q-1}2\cdot
\f d{2^{m+1}}}\Qs{qy^2}{\f{x+c}{2^m}+\f d{2^m}i}\Qs
{y^2}{\f{x+c}{2^m}+\f d{2^m}i}
\\&=(-1)^{\f{q^2-1}8+\f{q-1}2\cdot
\f d{2^{m+1}}}\Qs{c^2+d^2-x^2}{\f{x+c}{2^m}+\f d{2^m}i} \Ls
y{\f{(x+c)^2+d^2}{2^{2m}}}\endalign$$ and that
$$\align &\Qs{c^2+d^2-x^2}{\f{x+c}{2^m}+\f d{2^m}i}\\&
=\Qs{-2x(x+c)+d^2+(x+c)^2}{\f{x+c}{2^m}+\f d{2^m}i}
=\Qs{-2x(x+c)}{\f{x+c}{2^m}+\f d{2^m}i}
\\&=\Qs {2^{m+s+1}}{\f{x+c}{2^m}+\f d{2^m}i}
\Qs {-x_0}{\f{x+c}{2^m}+\f d{2^m}i} \Qs{(x+c)/2^m}{\f{x+c}{2^m}+\f
d{2^m}i}
\\&=i^{(-1)^{(\f{x+c}{2^m}-1)/2}
\f d{2^{m+1}}(m+s+1)}(-1)^{\f{x_0+1}2\cdot \f d{2^{m+1}}}
\Qs{\f{x+c}{2^m}+\f d{2^m}i}{x_0} \\&\q\times
(-1)^{\f{\f{x+c}{2^m}-1}2\cdot \f d{2^{m+1}}}\Qs {\f{x+c}{2^m}+\f
d{2^m}i}{\f{x+c}{2^m}}
\\&=i^{(-1)^{(\f{x+c}{2^m}-1)/2}
\f d{2^{m+1}}(m+s+1)} \cdot (-1)^{\f{\f{x+c}{2^m}+x_0}2\cdot \f
d{2^{m+1}}} \Qs{x+c+di}{x_0}\Qs {i}{\f{x+c}{2^m}}.\endalign$$ Thus,

$$\align &\Qs{c^2+d^2-x^2}{\f{x+c}{2^m}+\f d{2^m}i}
\\&=i^{(-1)^{(\f{x+c}{2^m}-1)/2} \f d{2^{m+1}}(m+s+1)}
(-1)^{\f{\f{x+c}{2^m}+x_0}2\cdot \f d{2^{m+1}}}\Ls
{2}{\f{x+c}{2^m}}\Qs{c+di}{x_0}
\\&=i^{(-1)^{(\f{x+c}{2^m}-1)/2}
\f d{2^{m+1}}(m+s+1)}   (-1)^{\f{\f{x+c}{2^m}+x_0}2\cdot \f
d{2^{m+1}}}\Ls {2}{\f{x+c}{2^m}}(-1)^{\f{x_0-1}2\cdot \f d2}\Qs
{x_0}{c+di}
\\&=i^{(-1)^{(\f{x+c}{2^m}-1)/2}
\f d{2^{m+1}}(m+s+1)}   (-1)^{\f{\f{x+c}{2^m}+x_0}2\cdot \f
d{2^{m+1}}}\Ls {2}{\f{x+c}{2^m}}\\&\q\times(-1)^{\f{x_0-1}2\cdot \f
d2}\Qs 2{c+di}^{-s}\Qs {x}{c+di}
\\&=(-1)^{\f{\f{x+c}{2^m}+x_0}2\cdot \f
d{2^{m+1}}+\f{x_0-1}2\cdot \f d2}i^{(-1)^{(\f{x+c}{2^m}-1)/2} \f
d{2^{m+1}}(m+s+1)-\f d2s}   \Ls {2}{\f{x+c}{2^m}}\Qs {x}{c+di}.
\endalign$$
By Lemma 2.11,
$$\Ls{y}{\f{(x+c)^2+d^2}{2^{2m}}}=(-1)^{\f{\f{(x+c)^2+d^2}{2^{2m}}-1}4t
}i^{\f d2t}\Qs{y^{-1}}{c+di}=(-1)^{\f d{2^{m+1}}t}i^{\f
d2t}\Qs{y^{-1}}{c+di}.$$ From the above we deduce
$$\aligned \Qs{d-(x+c)i}q&=(-1)^{\f{q^2-1}8+\f{q-1}2\cdot\f
d{2^{m+1}}}\cdot (-1)^{\f{\f{x+c}{2^m}+x_0}2\cdot \f
d{2^{m+1}}+\f{x_0-1}2\cdot \f d2}\\&\qq\times
i^{(-1)^{(\f{x+c}{2^m}-1)/2} \f d{2^{m+1}}(m+s+1)-\f d2s} \Ls
{2}{\f{x+c}{2^m}}
\\&\qq\times (-1)^{\f d{2^{m+1}}t}i^{\f
d2t}\Qs{x/y}{c+di}.\endaligned\tag 3.3$$ Clearly $m=0$ if and only
if $2\mid x$. If $2\mid x$, then $2\nmid (x+c)$, $m=0<r$, $2\nmid y$
and $q\e 1\mod 4$.  Thus, from (3.3) we deduce
$$\align \Qs{d-(x+c)i}q
&=(-1)^{\f{q^2-1}8}\cdot (-1)^{\f{x+c+x_0}2\cdot \f
d{2}+\f{x_0-1}2\cdot \f d2}\\&\qq\times i^{(-1)^{(x+c-1)/2} \f
d2(s+1)-\f d2s} \Ls {2}{x+c}\Qs{x/y}{c+di}
\\&=(-1)^{\f{q^2-1}8+(\f x2+1)\f d2}i^{\f d2}\Ls {2}{x+c}\Qs{x/y}{c+di}.
\endalign$$
Since $qy^2=d^2-(x+c)^2+2c(x+c)$ we see that $q\e qy^2\e
2(1-(-1)^{\f d2})-1+2c(x+c)\mod 8$. Thus, $x+c\e
\f{q+1}2-(1-(-1)^{\f d2})\mod 4$ and so $\f x2\e
\f{q-1}4-\f{1-(-1)^{\f d2}}2\mod 2$. We also have
 $$\align\Ls
2{x+c}&=(-1)^{\f{(x+c)^2-1}8}=(-1)^{\f{c^2-1}8}\cdot (-1)^{\f{\f
x2(\f x2+c)}2}\\&=(-1)^{\f{p-1-d^2}8}\cdot (-1)^{\f{\f x2(\f
x2+1)}2} =(-1)^{[\f p8]+[\f{x+2}4]}.\endalign$$ \par If $2\mid x$
and $p\e 1\mod 8$, then $4\mid d$ and so
$$\align \Qs{d-(x+c)i}q
&=(-1)^{\f{q^2-1}8+(\f x2+1)\f d2}i^{\f d2}\Ls
{2}{x+c}\Qs{x/y}{c+di}
\\&=(-1)^{\f{q-1}4+\f d4}\cdot
(-1)^{\f{p-1}8+[\f{x+2}4]}\Qs{x/y}{c+di}
\\&=(-1)^{\f x2+\f d4+\f{p-1}8+[\f{x+2}4]}\Qs{x/y}{c+di}
\\&=(-1)^{\f{p-1}8+\f d4+[\f x4]}\Qs{x/y}{c+di}.\endalign$$
Recall that $\sqs{d-(x+c)i}q=i^k$. We then get
$\sqs{x/y}{c+di}=(-1)^{\f{p-1}8+\f d4+[\f x4]}i^k.$ Now applying
Lemma 2.13 we obtain
$$q^{\f{p-1}8}\e (-1)^{\f d4+[\f x4]}\Ls dc^k\mod p.$$
\par
   If $2\mid x$ and
$p\e 5\mod 8$, then $2\ \Vert\ d$ and so
$$\align \Qs{d-(x+c)i}q
&=(-1)^{\f{q^2-1}8+(\f x2+1)\f d2}i^{\f d2}\Ls
{2}{x+c}\Qs{x/y}{c+di}
\\&=(-1)^{\f{q-1}4+\f x2+1}i\cdot
(-1)^{\f{p-5}8+[\f{x+2}4]}\Qs{x/y}{c+di}
\\&=(-1)^{\f{p-5}8+[\f{x+2}4]}i\Qs{x/y}{c+di}.\endalign$$
Recall that $\sqs{d-(x+c)i}q=i^k$. We then get
$\sqs{x/y}{c+di}=(-1)^{\f{p-5}8+[\f {x+2}4]}i^{k-1}.$ Now applying
Lemma 2.13 we obtain
$$q^{\f{p-5}8}\e (-1)^{[\f {x+2}4]}\Ls dc^{k-1}\f yx\mod p.$$

\par If $m=1<r$, then $x\e 1\mod 4$, $s=0$, $4\mid d$ and $p\e 1\mod 8$.
Since $p\e x^2\e 1\mod 8$ we have $8\mid qy^2$ and so $4\mid y$.
Thus,
$$\align 1-c\cdot\f{x+c}2&\e \Ls{x+c}2^2-c\cdot\f{x+c}2
=4\Big(\Ls d4^2-q\Ls y4^2\Big)\\& \e 4\Big(\f d4+\f y4\Big)=d+y\mod
8\endalign$$ and so $\f{x+c}2\e c-d-y\mod 8$. Therefore $\f{x+c}2\e
1\mod 4$ and
$$\align\Ls
2{\f{x+c}2}=\Ls 2{c-d-y}=(-1)^{\f{c-1}4+\f d4+\f y4}
=(-1)^{\f{c^2-1}8+\f d4+\f y4}=(-1)^{\f{p-1}8+\f d4+\f y4}.
\endalign$$
Hence $$\align &(-1)^{\f{q^2-1}8+\f{q-1}2\cdot\f d{2^{m+1}}}\cdot
(-1)^{\f{\f{x+c}{2^m}+x_0}2\cdot \f d{2^{m+1}}+\f{x_0-1}2\cdot \f
d2}\\&\qq\times i^{(-1)^{(\f{x+c}{2^m}-1)/2} \f d{2^{m+1}}(m+s+1)-\f
d2s} \Ls {2}{\f{x+c}{2^m}} \cdot (-1)^{\f d{2^{m+1}}t}i^{\f d2t}
\\&=(-1)^{\f{q^2-1}8+\f{q-1}2\cdot \f d4}\cdot (-1)^{\f{\f{x+c}2+x}2
\cdot \f d4}\cdot (-1)^{\f d4}\cdot (-1)^{\f{p-1}8+\f d4+\f y4}\cdot
(-1)^{\f d4t}\cdot (-1)^{\f d4t}
\\&=(-1)^{\f{p-1}8+\f{q^2-1}8+\f{q+1}2\cdot \f d4+\f y4}.\endalign$$
This together with (3.3) gives
$$\Qs{d-(x+c)i}q=(-1)^{\f{p-1}8+\f{q^2-1}8+\f{q+1}2\cdot \f d4+\f y4}
\Qs{x/y}{c+di}.$$ Since $\sqs{d-(x+c)i}q=i^k$, we then get $\sqs
{x/y}{c+di}=(-1)^{\f{p-1}8+\f{q^2-1}8+\f{q+1}2\cdot \f d4+\f
y4}i^k$. Now applying Lemma 2.13 we obtain
$$q^{\f{p-1}8}\e (-1)^{\f{q^2-1}8+\f{q+1}2\cdot \f d4+\f
y4}\Ls dc^k\mod p.$$

\par Now we assume $2\le m<r$. Then $x\e 3\mod 4$, $s=0$,
$8\mid d$ and so $p\e 1\mod 8$.
 Since $qy^2=d^2-(x+c)^2+2c(x+c)$ we
see that
$$q\f{y^2}{2^{m+1}}=2^{2r-m-1}d_0^2-2^{m-1}\Ls{x+c}{2^m}^2
+c\cdot\f{x+c}{2^m}.$$ As $m\ge 2$ and $2r\ge 2(m+1)\ge m+4$, we
must have $2^{m+1}\ \Vert\ y^2$, $2\mid m+1$, $t=\f{m+1}2$ and $q\e
-2^{m-1}+c\cdot \f{x+c}{2^m}\mod 8$.  Thus $m\ge 3$, $\f{x+c}{2^m}\e
c(2^{m-1}+q)\e c-1+2^{m-1}+q\mod 8$ and so $\f{x+c}{2^m}\e q\mod 4$.
Therefore, by (3.3) we have
$$\aligned &\Qs{d-(x+c)i}q\\&=(-1)^{\f{q^2-1}8+\f{q-1}2\cdot\f
d{2^{m+1}}}\cdot (-1)^{\f{\f{x+c}{2^m}+x}2\cdot \f d{2^{m+1}}}
i^{(-1)^{(\f{x+c}{2^m}-1)/2} \f d{2^{m+1}}(m+1)}\\&\qq\times \Ls
{2}{\f{x+c}{2^m}}
 (-1)^{\f d{2^{m+1}}t}i^{\f
d2t}\Qs{x/y}{c+di}
\\&= (-1)^{\f{q^2-1}8}\cdot (-1)^{\f d{2^{m+1}}\cdot\f{m+1}2}
 \Ls 2{c-1+2^{m-1}+q}\cdot (-1)^{\f d{2^{m+1}}t+\f d4t}\Qs{x/y}{c+di}
\\&= (-1)^{\f{c-1}4+2^{m-3}}\Qs{x/y}{c+di}.
\endaligned$$
Note that $(-1)^{\f{p-1}8}=(-1)^{\f{c^2-1}8} =(-1)^{\f{c-1}4}$,
$\sqs{d-(x+c)i}q=i^k$ and
$$\aligned(-1)^{2^{m-3}}&=\cases 1&\t{if $m>3$,}
\\-1&\t{if $m=3$}\endcases
=\cases 1&\t{if $t>2$,}
\\-1&\t{if $t=2$}\endcases
\\&=(-1)^{\f y4}.\endaligned$$ We get
$\sqs {x/y}{c+di}=(-1)^{\f{p-1}8+\f y4}i^k.$ Recall that $8\mid d$.
Applying Lemma 2.13
 we obtain
$$q^{\f{p-1}8}\e (-1)^{\f y4}
\Ls dc^k=(-1)^{\f{q+1}2\cdot \f d4+\f y4} \Ls dc^k\mod p.$$ This
completes the proof.

\pro{Theorem 3.4} Let $p$ be a prime of the form $4k+1$, $q\in\Bbb
Z$, $2\nmid q$ and $p\nmid q$. Suppose that $p=c^2+d^2=x^2+qy^2$
with $c,d,x,y\in\Bbb Z$, $2\nmid x$, $c\e 1\mod 4$, $d=2^rd_0$,
$y=2^ty_0$, $d_0\e y_0\e 1\mod 4$, $(d_0,x+c)=1$, $2^r\ \Vert\
(x+c)$ and $\sqs{d/(x+c)-i}q=i^k$.
\par $(\t{\rm i})$ If $q\e 1\mod 4$, then
$$q^{[p/8]}\e
\cases (-1)^{\f{q+1}2\cdot \f d4+\f y4}\sls dc^k=(-1)^{\f y4}\sls
dc^k\mod p&\t{if $p\e 1\mod 8$,}
\\(-1)^{\f{q-1}4}\sls dc^{k-1}\f yx\mod p&\t{if $p\e 5\mod 8$.}
\endcases$$
\par
$(\t{\rm ii})$ If $q\e 3\mod 4$, then
$$q^{[p/8]}\e
\cases (-1)^{\f{q+1}2\cdot \f d4+\f y4}\sls dc^k=(-1)^{\f y4}\sls
dc^k\mod p&\t{if $p\e 1\mod 8$,}
\\(-1)^{\f{q+1}4}\sls dc^k \f yx\mod p&\t{if $p\e 5\mod 8$.}
\endcases$$
  \endpro

Proof.   As $(x+c,d_0)=1$, by Lemma 2.12 we have $(qy_0,x+c)=1$ and
$(qy_0^2,(x+c)^2+d^2)=1$. Note that $(x,y)^2\mid p$. We also have
$(x,y)=1$. Since $qy^2=d^2-(x+c)^2+2c(x+c)$ we see that
$$q\f{y^2}{2^{r+1}}=2^{r-1}\Big(d_0^2-\Ls{x+c}{2^r}^2\Big)
+c\cdot\f{x+c}{2^r}.$$ As $2^r\ \Vert\ (x+c)$ and $r\ge 1$, we must
have $2^{r+1}\ \Vert\ y^2$, $2\mid (r+1)$, $t=\f{r+1}2$ and $q\e
qy_0^2 \e c\cdot \f{x+c}{2^r}\mod 8$. Thus $\f{x+c}{2^r}\e cq\e
c-1+q\mod 8$ and so $\f{x+c}{2^r}\e q\mod 4$.  We also have
$$\align\Qs{d-(x+c)i}q&=\Qs{-i}q\Qs{x+c+di}q=
\Ls 2q\Qs{\f{x+c}{2^r}+\f d{2^r}i}q \\&=\Ls 2q\Qs{1+i}q\Qs
{\f{x+c+d}{2^{r+1}}-\f{x+c-d}{2^{r+1}}i}q.\endalign$$ We first
assume $q\e 1\mod 4$. Then $\f{x+c}{2^r}\e q\e 1\e d_0=\f d{2^r}\mod
4$ and so $2^{r+2}\mid (x+c-d)$. Using Lemmas 2.1-2.5 we see that
$$\align&\Qs {\f{x+c+d}{2^{r+1}}-\f{x+c-d}{2^{r+1}}i}q
\\&=\Qs
q{\f{x+c+d}{2^{r+1}}-\f{x+c-d}{2^{r+1}}i} =\Qs
{qy^2}{\f{x+c+d}{2^{r+1}}-\f{x+c-d}{2^{r+1}}i} \Qs
{y^2}{\f{x+c+d}{2^{r+1}}-\f{x+c-d}{2^{r+1}}i}
\\&=\Qs {-2x(x+c)+\f
12((x+c+d)^2+(x+c-d)^2)}{\f{x+c+d}{2^{r+1}}-\f{x+c-d}{2^{r+1}}i}
 \Ls y{\f{(x+c+d)^2+(x+c-d)^2}{2^{2r+2}}}
\\&=\Qs{-2x(x+c)}{\f{x+c+d}{2^{r+1}}-\f{x+c-d}{2^{r+1}}i} \Ls
y{((x+c)^2+d^2)/2^{2r+1}}
\\&=
\Qs 2{\f{x+c+d}{2^{r+1}}-\f{x+c-d}{2^{r+1}}i}^{r+1} \\&\qq\times\Qs
{-x}{\f{x+c+d}{2^{r+1}}-\f{x+c-d}{2^{r+1}}i} \Qs{(x+c)/2^r}
{\f{x+c+d}{2^{r+1}}-\f{x+c-d}{2^{r+1}}i}\Ls
y{((x+c)^2+d^2)/2^{2r+1}}
\\&= i^{(-1)^{\f{\f{x+c+d}{2^{r+1}}-1}2}\f{x+c-d}{2^{r+2}}(r+1)}
\cdot (-1)^{\f{-x-1}2\cdot \f{x+c-d}{2^{r+2}}} \Qs
{\f{x+c+d}{2^{r+1}}-\f{x+c-d}{2^{r+1}}i}x
 \\&\qq\times
 \Qs {\f{x+c+d}{2^{r+1}}-\f{x+c-d}{2^{r+1}}i}{\f{x+c}{2^r}}
 \Ls y{((x+c)^2+d^2)/2^{2r+1}}
\\&=
(-1)^{\f{x+c-d}{2^{r+2}}( \f{r+1}2-\f{x+1}2)}
 \Qs {c+d-(c-d)i}x
 \\&\qq\times
 \Qs {d+di}{(x+c)/2^r}
 \Ls y{((x+c)^2+d^2)/2^{2r+1}}.
\endalign$$
As $\f{x+c}{2^r}\e 1\mod 4$, $c+d-(c-d)i=i^3(1+i)(c+di)$ and
$\f{x-1}2\cdot \f d2\e 0\mod 2$, we have
$$\align &\Qs {c+d-(c-d)i}x
   \Qs {d+di}{(x+c)/2^r}
\\&=\Qs ix^3\Qs{1+i}x\Qs{c+di}x\Qs{1+i}{(x+c)/2^r}
\\&=(-1)^{\f{x-\sls{-1}x}4}\cdot i^{\f{\sls{-1}xx-1}4}\cdot
\Qs x{c+di} \cdot i^{\f{\f{x+c}{2^r}-1}4}.\endalign$$ By Lemma 2.11,
$$\aligned \Ls y{((x+c)^2+d^2)/2^{2r+1}}&
=(-1)^{\f{\f{(x+c)^2+d^2}{2^{2r+1}}-1}4t+\f{y_0-1}4(2r+1)}i^{\f
d2t}\Qs{y^{-1}}{c+di}
\\&=(-1)^{\f{\f{x+c}{2^r}+\f d{2^r}-2}4
\cdot\f{r+1}2+\f{y_0-1}4}\cdot i^{\f d2\cdot \f{r+1}2}
\Qs{y^{-1}}{c+di}
\\&=(-1)^{\f{x+c-d}{2^{r+2}}
\cdot\f{r+1}2+\f{y_0-1}4}\cdot
i^{(1-(-1)^{\f{p-1}4})/2}\Qs{y^{-1}}{c+di}.
\endaligned$$
Thus, from the above we derive
$$\align\Qs{d-(x+c)i}q&=\Ls 2q\Qs{1+i}q\Qs
{\f{x+c+d}{2^{r+1}}-\f{x+c-d}{2^{r+1}}i}q
\\&=(-1)^{\f{q-1}4}i^{\f{q-1}4} \cdot (-1)^{\f{x+c-d}{2^{r+2}}\cdot
\f{r-x}2} \cdot (-1)^{\f{x-\sls{-1}x}4}\cdot
i^{\f{\sls{-1}xx-1}4+\f{\f{x+c}{2^r}-1}4}
\\&\qq\times (-1)^{\f{x+c-d}{2^{r+2}}
\cdot\f{r+1}2+\f{y_0-1}4}\cdot i^{(1-(-1)^{\f{p-1}4})/2}
\Qs{x/y}{c+di} .\endalign$$ If $p\e 1\mod 8$, then $r\ge 2$, $x\e
3\mod 4$ and $r=2t-1\ge 3$.
 As $qy_0^2\e c\cdot\f{x+c}{2^r}\mod{16}$
and $8\mid (x+c)$ we have
$$\align\f{\f{x+c}{2^r}-1}4&\e \f{\f qcy_0^2-1}4
=\f 1c\Big(\f {q-c}4+\f{y_0^2-1}4+\f{q-1}4(y_0^2-1)\Big)
\\&\e \f {q-c}4+\f{y_0^2-1}4\mod 4.\endalign$$ and so
$$\align \Qs{d-(x+c)i}q
&=(-1)^{\f{q-1}4+\f{x+1}4+\f{y_0-1}4} i^{\f{q-1}4-\f{x+1}4+\f
{q-c}4+\f{y_0^2-1}4} \Qs{x/y}{c+di}
\\&=(-1)^{\f{x+1}4}i^{-\f{x+c}4} \Qs{x/y}{c+di}
=(-1)^{\f{c-1}4}i^{-\f{x+c}4} \Qs{x/y}{c+di}
\\&=(-1)^{\f{c^2-1}8}i^{-2^{r-2}(q-c+y_0^2)}\Qs{x/y}{c+di}
=(-1)^{\f{p-1}8}\cdot (-1)^{2^{r-3}}\Qs{x/y}{c+di}
\\&=(-1)^{\f{p-1}8+\f y4}\Qs{x/y}{c+di}.
\endalign$$Since $\sqs{d-(x+c)i}q=i^k$, we get $\sqs{x/y}{c+di}=
(-1)^{\f{p-1}8+\f y4}i^k$. Now applying Lemma 2.13 we obtain
$q^{\f{p-1}8}\e (-1)^{\f y4}\sls dc^k\mod p.$
 If $p\e
5\mod 8$, then $r=1$, $x\e 1\mod 4$ and $t=\f{r+1}2=1$. From the
above we have
$$\align&\Qs{d-(x+c)i}q
\\&=(-1)^{\f{q-1}4}i^{\f{q-1}4} \cdot (-1)^{\f{x+c-d}{2^{r+2}}\cdot
\f{r-x}2} \cdot (-1)^{\f{x-\sls{-1}x}4}\cdot
i^{\f{\sls{-1}xx-1}4+\f{\f{x+c}{2^r}-1}4}
\\&\qq\times (-1)^{\f{x+c-d}{2^{r+2}}
\cdot\f{r+1}2+\f{y_0-1}4}\cdot i^{(1-(-1)^{\f{p-1}4})/2}
\Qs{x/y}{c+di}
\\&=(-1)^{\f{q-1}4}i^{\f{q-1}4} \cdot
(-1)^{\f{x-1}4}i^{\f{x-1}4+\f{\f{x+c}2-1}4} \cdot
(-1)^{\f{x+c-d}8+\f{y_0-1}4}i\Qs{x/y}{c+di} .\endalign$$ Since
$qy_0^2=d_0^2-\sls{x+c}2^2+c\cdot\f{x+c}2$ we see that
$$\align i^{\f{\f{x+c}2-1}4}&=i^{\f{c\cdot \f{x+c}2-c}4}
=i^{\f{qy_0^2-d_0^2+\sls{x+c}2^2-c}4}
=(-1)^{\f{y_0^2-1}8+\f{\sls{x+c}2^2-d_0^2}8} i^{\f{q-c}4}
\\&=(-1)^{\f{y_0-1}4+\f{\f{x+c}2-d_0}4}i^{\f{q-c}4}
=(-1)^{\f{y_0-1}4+\f{x+c-d}8}i^{\f{q-1}4-\f{c-1}4}.\endalign$$ Thus,
$$\align \Qs{d-(x+c)i}q&=(-1)^{\f{q-1}4}i^{\f{q-1}4} \cdot
(-1)^{\f{x-1}4}i^{\f{x-1}4}\cdot
i^{\f{q-1}4-\f{c-1}4+1}\Qs{x/y}{c+di}
\\&=(-1)^{\f{x-1}4}i^{1+\f{x-c}4}\Qs{x/y}{c+di}.\endalign$$
Since $\f{x+c}2\e c-1+q\mod 8$ we have $x-c\e 2(q-1)\mod{16}$ and so
$i^{\f{x-c}4}=i^{\f{2(q-1)}4}=(-1)^{\f{q-1}4}$. Thus,
$$\Qs{d-(x+c)i}q=(-1)^{\f{x-c}4}i^{\f{x-c}4}
\cdot
(-1)^{\f{c-1}4}i\Qs{x/y}{c+di}=(-1)^{\f{q-1}4+\f{c-1}4}i\Qs{x/y}{c+di}$$
and so
$$\Qs{x/y}{c+di}=(-1)^{\f{q-1}4+\f{c-1}4}i^{k-1}
=(-1)^{\f{q-1}4+\f{p-d^2-1}8}i^{k-1}=(-1)^{\f{p-5}8+\f{q-1}4}i^{k-1}.$$
By appealing to Lemma 2.13 we deduce $q^{\f{p-5}8}\e
(-1)^{\f{q-1}4}\sls dc^{k-1}\f yx\mod p.$

\par Now we assume  $q\e 3\mod 4$. Then $\f{x+c}{2^r}\e q\e -1\e -d_0=-\f d{2^r}\mod
4$ and so $2^{r+2}\mid (x+c+d)$. Hence, using Lemmas 2.1-2.5 we get
$$\align&\Qs {\f{x+c+d}{2^{r+1}}-\f{x+c-d}{2^{r+1}}i}q
\\&=\Qs{-i}q\Qs {\f{x+c-d}{2^{r+1}}+\f{x+c+d}{2^{r+1}}i}q
=(-1)^{\f{q+1}4}\Qs {-q}{\f{x+c-d}{2^{r+1}}+\f{x+c+d}{2^{r+1}}i}
\\&=(-1)^{\f{q+1}4}\Qs {-qy^2}{\f{x+c-d}{2^{r+1}}+\f{x+c+d}{2^{r+1}}i}
\Qs {y^2}{\f{x+c-d}{2^{r+1}}+\f{x+c+d}{2^{r+1}}i}
\\&=(-1)^{\f{q+1}4}\Qs {2x(x+c)-\f
12((x+c+d)^2+(x+c-d)^2)}{\f{x+c-d}{2^{r+1}}+\f{x+c+d}{2^{r+1}}i}
\\&\qq\times \Ls y{((x+c-d)^2+(x+c+d)^2)/2^{2r+2}}
\\&=(-1)^{\f{q+1}4}\Qs{2^{r+1}(-x)(-\f{x+c}{2^r})}{\f{x+c-d}{2^{r+1}}+\f{x+c+d}{2^{r+1}}i} \Ls
y{((x+c)^2+d^2)/2^{2r+1}}\endalign$$ and therefore
$$\align &\Qs {\f{x+c+d}{2^{r+1}}-\f{x+c-d}{2^{r+1}}i}q
\\&=(-1)^{\f{q+1}4} \Qs
2{\f{x+c-d}{2^{r+1}}+\f{x+c+d}{2^{r+1}}i}^{r+1}\Qs
{-x}{\f{x+c-d}{2^{r+1}}+\f{x+c+d}{2^{r+1}}i}\\&\qq\times
 \Qs{-(x+c)/2^r}
{\f{x+c-d}{2^{r+1}}+\f{x+c+d}{2^{r+1}}i}\Ls
y{((x+c)^2+d^2)/2^{2r+1}}
\\&= (-1)^{\f{q+1}4}\cdot
i^{(-1)^{\f{\f{x+c-d}{2^{r+1}}-1}2}\f{x+c+d}{2^{r+2}}(r+1)}
 (-1)^{\f{-x-1}2\cdot \f{x+c+d}{2^{r+2}}}
 \Qs {\f{x+c-d}{2^{r+1}}+\f{x+c+d}{2^{r+1}}i}x
 \\&\qq\times
 \Qs {\f{x+c-d}{2^{r+1}}+\f{x+c+d}{2^{r+1}}i}{\f{x+c}{2^r}}
 \Ls y{((x+c)^2+d^2)/2^{2r+1}}
\\&=(-1)^{\f{q+1}4+\f{x+c+d}{2^{r+2}}(\f{r+1}2- \f{x+1}2)}
 \Qs {c-d+(c+d)i}x
 \\&\qq\times
 \Qs {-d+di}{(x+c)/2^r}
 \Ls y{((x+c)^2+d^2)/2^{2r+1}}.
\endalign$$
As $\f{x+c}{2^r}\e 3\mod 4$, $\f{x-1}2\cdot \f d2\e 0\mod 2$,
$c-d+(c+d)i=(1+i)(c+di)$ and $-d+di=di(1+i)$, we have
$$\align &\Qs {c-d+(c+d)i}x
   \Qs {-d+di}{(x+c)/2^r}
\\&=\Qs{1+i}x\Qs{c+di}x\Qs i{(x+c)/2^r}\Qs{1+i}{(x+c)/2^r}
\\&=i^{\f{\sls {-1}xx-1}4} \Qs x{c+di}\cdot (-1)^{\f{\f{x+c}{2^r}+1}4}
\cdot i^{-\f{\f{x+c}{2^r}+1}4}.\endalign$$ By Lemma 2.11,
$$\align \Ls y{((x+c)^2+d^2)/2^{2r+1}}&
=(-1)^{\f{\f{(x+c)^2+d^2}{2^{2r+1}}-1}4t+\f{y_0-1}4(2r+1)}i^{\f
d2t}\Qs{y^{-1}}{c+di} \\&=(-1)^{\f{\f{x+c}{2^r}-\f d{2^r}+2}4
\cdot\f{r+1}2+\f{y_0-1}4}i^{\f d2\cdot \f{r+1}2}\Qs{y^{-1}}{c+di}
\\&=(-1)^{\f{x+c+d}{2^{r+2}}\cdot\f{r+1}2+\f{y_0-1}4}i^{(1-(-1)^{\f{p-1}4})/2
}\Qs{y^{-1}}{c+di} .
\endalign$$
Thus, from the above we derive
$$\align&\Qs{d-(x+c)i}q\\&=\Ls 2q\Qs{1+i}q\Qs
{\f{x+c+d}{2^{r+1}}-\f{x+c-d}{2^{r+1}}i}q
\\&=(-1)^{\f{q+1}4} i^{-\f{q+1}4}
\cdot (-1)^{\f{q+1}4+\f{x+c+d}{2^{r+2}}( \f{r+1}2-\f{x+1}2)} \cdot
(-1)^{\f{\f{x+c}{2^r}+1}4}i^{\f{\sls{-1}xx-1}4-\f{\f{x+c}{2^r}+1}4}
\\&\qq\times
(-1)^{\f{x+c+d}{2^{r+2}}\cdot\f{r+1}2+\f{y_0-1}4}i^{(1-(-1)^{\f{p-1}4})/2}
\Qs{x/y}{c+di}
\\&=(-1)^{\f{x+c+d}{2^{r+2}}\cdot \f{x+1}2
+\f{\f{x+c}{2^r}+1}4+\f{y_0-1}4}
i^{\f{1-(-1)^{(p-1)/4}}2-\f{q+1}4+\f{\sls{-1}xx-1}4-\f{\f{x+c}{2^r}+1}4}
\Qs{x/y}{c+di}.\endalign$$ For $p\e 1\mod 8$ we have $4\mid d$, $x\e
3\mod 4$ and $r=2t-1\ge 3$. As $qy_0^2\e c\cdot\f{x+c}{2^r}\mod{16}$
we have
$$\align\f{\f{x+c}{2^r}+1}4&\e \f{\f qcy_0^2+1}4
=\f 1c\Big(\f {c+q}4-\f{y_0^2-1}4+\f{q+1}4(y_0^2-1)\Big)
\\&\e \f {c+q}4-\f{y_0^2-1}4\mod 4
\endalign$$
and so

$$\align \Qs{d-(x+c)i}q
&=(-1)^{\f{\f{x+c}{2^r}+1}4+\f{y_0-1}4}
i^{-\f{q+1}4-\f{x+1}4-\f{\f{x+c}{2^r}+1}4} \Qs{x/y}{c+di}
\\&=(-1)^{\f{c+q}4+\f{y_0-1}4}
i^{-\f{q+1}4-\f{x+1}4-\f{c+q}4+\f{y_0^2-1}4} \Qs{x/y}{c+di}
\\&=(-1)^{\f{c-1}4}i^{-\f{x+c}4}\Qs{x/y}{c+di}
.\endalign$$ Clearly $(-1)^{\f{c-1}4}=(-1)^{\f{c^2-1}8}
=(-1)^{\f{p-1}8}$. Since $r>3\iff t=\f{r+1}2>2\iff 8\mid y$, we have
$i^{-\f{x+c}4}=i^{-2\cdot 2^{r-3}\cdot \f{x+c}{2^r}}
=(-1)^{2^{r-3}}=(-1)^{\f y4}.$ Note that $\sqs{d-(x+c)i}q=i^k$. We
then get $\sqs{x/y}{c+di}=(-1)^{\f{p-1}8+\f y4}i^k$. Now applying
Lemma 2.13 we deduce $q^{\f{p-1}8}\e (-1)^{\f y4}\sls dc^k\mod p$.
\par For $p\e 5\mod 8$ we have $r=1$, $x\e 1\mod 4$ and
$t=\f{r+1}2=1$. Thus,
$$\Qs{d-(x+c)i}q
=(-1)^{\f{x+c+d}8+\f{x+c+2}8+\f{y_0-1}4}i^{1-\f{q+1}4+\f{x-1}4-\f{x+c+2}8}
\Qs{x/y}{c+di}.$$ Since $qy_0^2=d_0^2-\sls{x+c}2^2+c\cdot\f{x+c}2$
we see that
$$\align i^{-\f{x+c+2}8}&=i^{-\f{c\cdot \f{x+c}2+c}4}
=i^{-\f{qy_0^2-d_0^2+\sls{x+c}2^2+c}4}
=(-1)^{\f{y_0^2-1}8+\f{\sls{x+c}2^2-d_0^2}8} i^{-\f{q+c}4}
\\&=(-1)^{\f{y_0-1}4+\f{\f{x+c}2+d_0}4}i^{-\f{q+c}4}
=(-1)^{\f{y_0-1}4+\f{x+c+d}8}i^{-\f{q+1}4-\f{c-1}4}.\endalign$$
Hence
$$\Qs{d-(x+c)i}q
=(-1)^{\f{q+1}4+\f{x+c+2}8}i^{1+\f{x-c}4}\Qs{x/y}{c+di}.$$ Since
$\f{x+c}2\e c-1+q\mod 8$ we have $x-c\e 2(q-1)\mod{16}$ and so
$i^{1+\f{x-c}4}=i^{1+\f{2(q-1)}4}=(-1)^{\f{q+1}4}$. We also have
$(-1)^{\f{x+c+2}8}=(-1)^{\f{2(q-1)+2c+2}8} =(-1)^{\f{q+1}4+\f{c-1}4}
=(-1)^{\f{q+1}4+\f{c^2-1}8}=(-1)^{\f{q+1}4+\f{p-5}8}$.
 Thus,
$$\Qs{d-(x+c)i}q=(-1)^{\f{x+c+2}8}\Qs{x/y}{c+di}
=(-1)^{\f{q+1}4+\f{p-5}8}\Qs{x/y}{c+di}$$ and so
$\sqs{x/y}{c+di}=(-1)^{\f{q+1}4}i^k.$ Now applying Lemma 2.13 we get
$q^{\f{p-5}8}\e (-1)^{\f{q+1}4}$ $ \sls dc^k\f yx\mod p.$ This
completes the proof.
 \pro{Theorem 3.5} Let $p$ be a prime of the
form $4k+1$, $q\in\Bbb Z$, $2\nmid q$ and $p\nmid q$. Suppose that
$p=c^2+d^2=x^2+qy^2$ with $c,d,x,y\in\Bbb Z$, $x\e 3\mod 4$, $c\e
1\mod 4$, $d=2^rd_0$, $y=2^ty_0$, $d_0\e y_0\e 1\mod 4$,
$(d_0,x+c)=1$ and $\sqs{d/(x+c)-i}q=i^k$.
\par $(\t{\rm i})$ If $p\e 1\mod 8$ and $2^{r+1}\mid (x+c)$, then
$$q^{\f{p-1}8}\e
 (-1)^{\f {q+1}2\cdot \f d4+\f y4}\Ls dc^k\mod p.$$
\par
$(\t{\rm ii})$ If $p\e 5\mod 8$, then
$$q^{\f{p-5}8}\e
\cases -\sls dc^{k-1}\f yx\mod p&\t{if $q\e 1\mod 4$,}
\\-\sls dc^k\f yx\mod p&\t{if $q\e 3\mod 4$.}
\endcases$$
  \endpro
Proof. Suppose $2^m\ \Vert\ (x+c)$. We may assume $m>r$. Using
Lemmas 2.1-2.5 we see that
$$\align \Qs{d-(x+c)i}q&=
\Qs{\f d{2^r}-\f{x+c}{2^r}i}q=(-1)^{\f{q-1}2\cdot\f{x+c}{2^{r+1}}}
\Qs q {\f d{2^r}-\f{x+c}{2^r}i}
\\&=(-1)^{\f{q-1}2\cdot\f{x+c}{2^{r+1}}}
\Qs {qy^2} {\f d{2^r}-\f{x+c}{2^r}i} \Qs {y^2} {\f
d{2^r}-\f{x+c}{2^r}i}
\\&=(-1)^{\f{q-1}2\cdot\f{x+c}{2^{r+1}}}
\Qs {(x+c)^2+d^2-2x(x+c)} {\f d{2^r}-\f{x+c}{2^r}i} \Ls y{\sls
d{2^r}^2+\sls{x+c}{2^r}^2}
\\&=(-1)^{\f{q-1}2\cdot 2^{m-r-1}}
\Qs {-2x(x+c)} {\f d{2^r}-\f{x+c}{2^r}i} \Ls
y{((x+c)^2+d^2)/2^{2r}}\endalign$$ and
$$\align \Qs {-2x(x+c)} {\f d{2^r}-\f{x+c}{2^r}i}
&=\Qs 2 {\f d{2^r}-\f{x+c}{2^r}i}^{m+1} \Qs {-x} {\f
d{2^r}-\f{x+c}{2^r}i} \Qs {(x+c)/2^m} {\f d{2^r}-\f{x+c}{2^r}i}
\\&=i^{-\f{x+c}{2^{r+1}}(m+1)}\Qs {\f d{2^r}-\f{x+c}{2^r}i}
x\cdot (-1)^{\f{\f{x+c}{2^m}-1}2\cdot \f{x+c}{2^{r+1}}} \Qs {\f
d{2^r}-\f{x+c}{2^r}i}{(x+c)/2^m}
\\&=(-1)^{\f{\f{x+c}{2^m}-1}2\cdot \f{x+c}{2^{r+1}}}
i^{-\f{x+c}{2^{r+1}}(m+1)}\Qs{d-(x+c)i}x\Qs{d/2^r}{(x+c)/2^m}
\\&=(-1)^{\f{\f{x+c}{2^m}-1}2\cdot \f{x+c}{2^{r+1}}}
i^{-\f{x+c}{2^{r+1}}(m+1)}\Qs{d-ci}x
\\&=(-1)^{\f{\f{x+c}{2^m}-1}2\cdot \f{x+c}{2^{r+1}}}
i^{-\f{x+c}{2^{r+1}}(m+1)}\Qs{-i}x\Qs{c+di}x
\\&=(-1)^{\f{\f{x+c}{2^m}-1}2\cdot \f{x+c}{2^{r+1}}+\f{x+1}4+\f d2}
i^{-\f{x+c}{2^{r+1}}(m+1)}\Qs x{c+di}.\endalign$$
By Lemma 2.11,
$$\align\Ls y{((x+c)^2+d^2)/2^{2r}}
=(-1)^{\f{\f{(x+c)^2+d^2}{2^{2r}}-1}4t} i^{\f d2t}\Qs{y^{-1}}{c+di}
=(-1)^{\f{x+c}{2^{r+1}}t}i^{\f d2t} \Qs{y^{-1}}{c+di}.\endalign$$
Therefore,
$$\align&\Qs{d-(x+c)i}q\\&=(-1)^{\f {q-1}2\cdot 2^{m-r-1}}\cdot
(-1)^{\f{\f{x+c}{2^m}-1}2\cdot \f{x+c}{2^{r+1}}+\f{x+1}4+\f d2}
i^{-\f{x+c}{2^{r+1}}(m+1)}(-1)^{\f{x+c}{2^{r+1}}t}i^{\f d2t}
\Qs{x/y}{c+di}.\endalign$$ Since $qy^2=2c(x+c)-(x+c)^2+d^2$ we see
that
$$q\f{y^2}{2^{m+1}}=c\cdot\f{x+c}{2^m}-2^{m-1}\Ls{x+c}{2^m}^2
+2^{2r-m-1}d_0^2.$$ If $m<2r-1$, then $r>1$, $4\mid d$ and $p\e
1\mod 8$. As $2^m\ \Vert\ (x+c)$ and $m\ge 2$ we have $2^{m+1}\
\Vert\ y^2$, $2\mid (m+1)$, $t=\f{m+1}2$ and $q\e qy_0^2\e c\cdot
\f{x+c}{2^m}-2^{m-1}+2^{2r-m-1}\mod 8$. Thus $m\ge 3$,
$\f{x+c}{2^m}\e c(2^{m-1}+q-2^{2r-m-1})\e q-2^{2r-m-1}\mod 4$. For
$m=r+1<2r-1$ we have $2\mid r$, $r>\f{m+1}2\ge 2$, $r\ge 4$ and
$t=\f{r+2}2\ge 3$. Thus,
$$\align &(-1)^{\f {q-1}2\cdot 2^{m-r-1}}\cdot
(-1)^{\f{\f{x+c}{2^m}-1}2\cdot \f{x+c}{2^{r+1}}+\f{x+1}4+\f d2}
i^{-\f{x+c}{2^{r+1}}(m+1)}(-1)^{\f{x+c}{2^{r+1}}t}i^{\f d2t}
\\&=(-1)^{\f{q-1}2}\cdot (-1)^{\f{q-1-2^{r-2}}2+\f{x+1}4}
i^{(2^{r-2}-q)(r+2)}(-1)^{\f{r+2}2(1+\f d4)}
\\&=(-1)^{2^{r-3}+\f{x+1}4}\cdot(-1)^{(2^{r-2}-q)\f{r+2}2}
\cdot(-1)^{\f{r+2}2}
\\&=(-1)^{2^{r-3}+\f{x+1}4}=(-1)^{\f{x+1}4}=(-1)^{\f{1-c}4}
=(-1)^{\f{c^2-1}8}
\\&=(-1)^{\f{p-1}8}=(-1)^{\f{p-1}8+\f{q+1}2\cdot \f d4+\f y4}.
\endalign$$
For $r+2\le m<2r-1$ we have $r\ge 4$, $m\ge r+2\ge 6$,
$t=\f{m+1}2\ge 4$ and so
$$\align &(-1)^{\f {q-1}2\cdot 2^{m-r-1}}\cdot
(-1)^{\f{\f{x+c}{2^m}-1}2\cdot \f{x+c}{2^{r+1}}+\f{x+1}4+\f d2}
i^{-\f{x+c}{2^{r+1}}(m+1)}(-1)^{\f{x+c}{2^{r+1}}t}i^{\f d2t}
\\&= (-1)^{\f{x+1}4}\cdot (-1)^{\f{x+c}{2^{r+2}}(m+1)}
\cdot (-1)^{\f d4t}=(-1)^{\f{x+1}4}
\\&=(-1)^{\f{1-c}4}=(-1)^{\f{c^2-1}8}=(-1)^{\f{p-1}8}
=(-1)^{\f{p-1}8+\f{q+1}2\cdot \f d4+\f y4}.
\endalign$$
Now we assume $m\ge 2r-1$.  As $qy^2=2c(x+c)-(x+c)^2+d^2$ we have
$$q\f{y^2}{2^{2r}}=c\cdot
\f{x+c}{2^{2r-1}}-2^{2r-2}\Ls{x+c}{2^{2r-1}}^2+d_0^2.$$ Thus,
$2^r\mid y$ and so $t\ge r$. For $r\ge 2$ we have $m\ge 3$ and so
$8\mid (x+c)$. If $m=2r-1$, then $2\mid \f {y^2}{2^{2r}}$ and so
$2^{r+1}\mid y$. Thus $t\ge r+1\ge 3$ and therefore
$\f{x+c}{2^{2r-1}}\e -\f{d_0^2}c\e 3\mod 4$. For $r=2$ and $m=3$, we
have
$$\align &(-1)^{\f {q-1}2\cdot 2^{m-r-1}}\cdot
(-1)^{\f{\f{x+c}{2^m}-1}2\cdot \f{x+c}{2^{r+1}}+\f{x+1}4+\f d2}
i^{-\f{x+c}{2^{r+1}}(m+1)}(-1)^{\f{x+c}{2^{r+1}}t}i^{\f d2t}
\\&= (-1)^{\f{q-1}2}\cdot (-1)^{1+\f{x+1}4}\cdot
(-1)^{(1+\f d4)t}=(-1)^{\f{q+1}2+\f{1-c}4}
\\&=(-1)^{\f{q+1}2+\f{c^2-1}8}=(-1)^{\f{p-1}8+\f{q+1}2}
=(-1)^{\f{p-1}8+\f{q+1}2\cdot\f d4+\f y4}.
\endalign$$
For $r>2$ and $m=2r-1$, we have $m>r+1$ and so
$$\align &(-1)^{\f {q-1}2\cdot 2^{m-r-1}}\cdot
(-1)^{\f{\f{x+c}{2^m}-1}2\cdot \f{x+c}{2^{r+1}}+\f{x+1}4+\f d2}
i^{-\f{x+c}{2^{r+1}}(m+1)}(-1)^{\f{x+c}{2^{r+1}}t}i^{\f d2t}
\\&=  (-1)^{\f{x+1}4}=(-1)^{\f{1-c}4}
=(-1)^{\f{c^2-1}8}=(-1)^{\f{p-1}8} =(-1)^{\f{p-1}8+\f{q+1}2\cdot\f
d4+\f y4}.
\endalign$$
For $m\ge 2r$, as $$q\f{y^2}{2^{2r}}=2^{m-2r+1}c\cdot \f{x+c}{2^m}
-2^{2m-2r} \Ls{x+c}{2^m}^2+d_0^2$$ we have $t=r$ and $q\e q\sls
y{2^r}^2\e 2^{m-2r+1}c\cdot\f{x+c}{2^m}+1\mod 4$. Thus,
$$q\e \cases 3\mod 4&\t{if $m=2r$,}
\\1\mod 4&\t{if $m>2r$}.\endcases$$ For $p\e 1\mod 8$ we have $r\ge
2$, $m\ge 2r>r+1$ and therefore
$$\aligned &(-1)^{\f {q-1}2\cdot 2^{m-r-1}}\cdot
(-1)^{\f{\f{x+c}{2^m}-1}2\cdot \f{x+c}{2^{r+1}}+\f{x+1}4+\f d2}
i^{-\f{x+c}{2^{r+1}}(m+1)}(-1)^{\f{x+c}{2^{r+1}}t}i^{\f d2t}
\\&=  (-1)^{\f{x+1}4}\cdot (-1)^{\f{x+c}{2^{r+2}}(m+1)}\cdot
(-1)^{\f d4r} =(-1)^{\f{1-c}4+2^{m-2-r}(m+1)}
\\&=(-1)^{\f{c^2-1}8+2^{m-2-r}(m+1)}=(-1)^{\f{p-1}8+2^{m-2-r}(m+1)}
\\&=\cases -(-1)^{\f{p-1}8} =(-1)^{\f{p-1}8+\f{q+1}2\cdot
\f d4+\f y4}&\t{if $r=2$ and
$m=4$,}\\(-1)^{\f{p-1}8}=(-1)^{\f{p-1}8+\f{q+1}2\cdot \f d4+\f
y4}&\t{if $r=2$ and $m>4$,}
\\ (-1)^{\f{p-1}8}=(-1)^{\f{p-1}8+\f{q+1}2\cdot
\f d4+\f y4}&\t{if $r>2$.}
\endcases
\endaligned$$
From the above we deduce that for $p\e 1\mod 8$, we always have
$$\Qs{d-(x+c)i}q=(-1)^{\f{p-1}8+\f{q+1}2\cdot
\f d4+\f y4}\Qs{x/y}{c+di}.$$ Since $\sqs{d-(x+c)i}q=i^k$, by the
above and Lemma 2.13 we derive
$$q^{\f{p-1}8}\e (-1)^{\f{q+1}2\cdot
\f d4+\f y4}\Ls dc^k\mod p.$$
 \par From now on we assume $p\e 5\mod 8$. Then $r=1$, and
$q\e 1$ or $3\mod 4$ according as $m>2$ or $m=2$.
 Therefore
$$\aligned &(-1)^{\f {q-1}2\cdot 2^{m-r-1}}\cdot
(-1)^{\f{\f{x+c}{2^m}-1}2\cdot \f{x+c}{2^{r+1}}+\f{x+1}4+\f d2}
i^{-\f{x+c}{2^{r+1}}(m+1)}(-1)^{\f{x+c}{2^{r+1}}t}i^{\f d2t}
\\&=(-1)^{\f {q-1}2\cdot 2^{m-2}}\cdot
(-1)^{\f{\f{x+c}{2^m}-1}2\cdot
 \f{x+c}4+\f{x+1}4+1}
i^{-\f{x+c}4(m+1)}(-1)^{\f{x+c}4}i
\\&=\cases (-1)^{\f{q-1}2+\f{x+c-4}8+\f{x+1}4+1+\f{x+c}4}i^{\f{x+c}4+1}
&\t{if $m=2$,}
\\(-1)^{\f{x+1}4+1+\f{x+c}8(m+1)+\f{x+c}4}i=-(-1)^{\f{c-1}4}i&\t{if $m>2$.}
\endcases
\endaligned$$
For $m=2$ we have $t=1$, $q\e 3\mod 4$,
$$q\e qy_0^2=2c\cdot \f{x+c}4-4\ls{x+c}4^2+d_0^2\e 2\cdot
\f{x+c}4-4+1\mod 8$$ and so $\f{x+c}4\e \f{q+3}2\mod 4$. Hence
$$\aligned &(-1)^{\f {q-1}2\cdot 2^{m-r-1}}\cdot
(-1)^{\f{\f{x+c}{2^m}-1}2\cdot \f{x+c}{2^{r+1}}+\f{x+1}4+\f d2}
i^{-\f{x+c}{2^{r+1}}(m+1)}(-1)^{\f{x+c}{2^{r+1}}t}i^{\f d2t}
\\&=
(-1)^{\f{q-1}2+\f{x+c-4}8+\f{x+1}4+1+\f{x+c}4}i^{\f{x+c}4+1}
=-(-1)^{\f{c-1}4}.\endaligned$$ Thus,
$$\Qs{d-(x+c)i}q=\cases -(-1)^{\f{c-1}4}i\sqs{x/y}{c+di}&\t{if $q\e
1\mod 4$,}
\\-(-1)^{\f{c-1}4}\sqs{x/y}{c+di}&\t{if $q\e
3\mod 4$.}\endcases$$ Since $\sqs{d-(x+c)i}q=i^k$ and
$(-1)^{\f{c-1}4}=(-1)^{\f{c^2-1}8}=(-1)^{\f{p-1-d^2}8}=(-1)^{\f{p-5}8}$,
by the above we obtain
$$\Qs{x/y}{c+di}=\cases -(-1)^{\f{p-5}8}i^{k-1}&\t{if $q\e 1\mod
4$,}\\-(-1)^{\f{p-5}8}i^k&\t{if $q\e 3\mod 4$.}
\endcases$$
Now applying Lemma 2.13 we obtain the result in the case $p\e 5\mod
8$.  Summarizing all the above we prove the theorem.
\par\q
\newline{\bf Remark 3.2}
We note that the $k$ in Theorems 3.3-3.5 depends only on $\f
d{x+c}\mod q$, and for $p\e 1\mod 8$ and $x\e 3\mod 4$ we always
have $$q^{\f{p-1}8}\e (-1)^{\f{q+1}2\cdot \f d4+\f y4}\Ls dc^k\mod
p.$$

\subheading{4. New reciprocity laws for quartic and octic residues}
\pro{Theorem 4.1} Let $p$ and $q$ be primes such that $p\e 1\mod 4$
and $q\e 3\mod 4$. Suppose $p=c^2+d^2=x^2+qy^2$, $c,d,x,y\in\Bbb Z$,
$c\e 1\mod 4$, $d=2^rd_0$, $y=2^ty_0$, $d_0\e y_0\e 1\mod 4$  and
$\sls{c-di}x^{\f{q+1}4}\e i^m\mod q$. Assume $(c,x+d)=1$ or
$(d_0,x+c)=1$. Then
$$q^{[p/8]}\e\cases (-1)^{\f y4+\f{q+1}4\cdot
\f{x-1}2}\sls dc^m\mod p&\t{if $p\e 1\mod
8$,}\\(-1)^{\f{q-3}4\cdot \f{x-1}2}\sls dc^m\f yx\mod p&\t{if $p\e
5\mod 8$.}
\endcases$$
\endpro
Proof. Since $p\e 1\mod 4$ and $q\e 3\mod 4$ we see that $q\nmid x$
and $x$ is odd. We first assume $(c,x+d)=1$. By Lemma 2.12 we have
$(q,(x+d)(c^2+(x+d)^2))=1$. It is easily seen that
$\f{c/(x+d)-i}{c/(x+d)+i}=\f{c-(x+d)i}{c+(x+d)i}\e \f{c-di}{ix}\mod
q$.
 Thus, for $k=0,1,2,3$, using Lemma 2.7 we get
$$\align &\Qs{c+(x+d)i}q=i^k\\&\iff\f c{x+d}\in Q_k(q)
\iff \Ls{\f c{x+d}-i}{\f c{x+d}+i}^{\f{q+1}4}\e i^k\mod q
\\&\iff
\Ls{c-di}{ix}^{\f{q+1}4}\e i^k\mod q \iff \Ls{c-di}x^{\f{q+1}4}\e
i^{\f{q+1}4+k}\mod q.
\endalign$$
Since  $\sls{c-di}x^{\f{q+1}4}\e i^m\mod q$, from the above we
deduce $$\Qs{c+(x+d)i}q=i^{m-\f{q+1}4}=\cases
(-1)^{\f{q+5}8}i^{m+1}&\t{if $q\e 3\mod
8$,}\\(-1)^{\f{q+1}8}i^m&\t{if $q\e 7\mod 8$.}\endcases$$ Now,
applying Theorem 3.2 we derive the result.
\par Now we assume $(d_0,x+c)=1$.
Since $qy^2=(x+c)^2-2x(x+c)+d^2$ we see that $(x+c,qy^2)=(x+c,d^2)$.
As $(d_0,x+c)=1$ we have $(q,x+c)=1$ and so
$(q,d^2+(x+c)^2)=(q,qy^2+2x(x+c))=(q,2x(x+c))=1$.
 It is easily seen
that $\f{d+(x+c)i}{d-(x+c)i}\e \f{c-di}{-x}\mod q$.
 Thus, for $k=0,1,2,3$, using Lemma 2.7 we get
$$\align &\Qs{d-(x+c)i}q=i^k\\&\iff-\f d{x+c}\in Q_k(q)
\iff \Ls{-\f d{x+c}-i}{-\f d{x+c}+i}^{\f{q+1}4}\e i^k\mod q
\\&\iff \Ls{d+(x+c)i}{d-(x+c)i}^{\f{q+1}4}\e i^k\mod q\iff
\Ls{c-di}{-x}^{\f{q+1}4}\e i^k\mod q
\\&\iff \Ls{c-di}x^{\f{q+1}4}\e i^{\f{q+1}2+k}\mod q.
\endalign$$
Since  $\sls{c-di}x^{\f{q+1}4}\e i^m\mod q$, from the above we
deduce that $\sqs{d-(x+c)i}q=i^{m-\f{q+1}2}=(-1)^{\f{q+1}4}i^m$. Now
applying Theorems 3.3-3.5 we deduce the result. The proof is now
complete.

\pro{Corollary 4.1} Let $p\e 1\mod 4$ and $q\e 3\mod{8}$ be primes
such that $p=c^2+d^2=x^2+qy^2$ with $c,d,x,y\in\Bbb Z$ and $q\mid
cd$. Suppose $c\e 1\mod 4$, $d=2^rd_0$, $y=2^ty_0$ and $d_0\e y_0\e
1\mod 4$. Assume $(c,d+x)=1$ or $(d_0,x+c)=1$.
\par $(\t{\rm i})$ If $p\e 1\mod 8$, then
$$q^{\f{p-1}8}\e \cases \pm (-1)^{\f{x-1}2+\f y4}\mod p&\t{if $x\e \pm c\mod{q}$,}
\\\mp(-1)^{\f{q-3}8+\f{x-1}2+\f y4}\f dc\mod p&\t{if $x\e\pm
d\mod{q}$.}
\endcases$$
\par $(\t{\rm ii})$ If $p\e 5\mod 8$, then
$$q^{\f{p-5}8}\e \cases \pm \f yx\mod p&\t{if $x\e\pm
c\mod{q}$,}\\\mp (-1)^{\f{q-3}8}\f{dy}{cx}\mod p&\t{if $x\e \pm
d\mod{q}$.}
\endcases$$
\endpro
Proof. If $x\e \pm c\mod q$, then $q\mid d$ and so
$\sls{c-di}x^{\f{q+1}4}\e (\pm 1)^{\f{q+1}4}= \pm 1\mod q$. If $x\e
\pm d\mod q$, then $q\mid c$ and so $\sls{c-di}x^{\f{q+1}4}\e (\mp
i)^{\f{q+1}4}= \mp (-1)^{\f{q-3}8}i\mod q$. Now applying Theorem 4.1
we deduce the result.
\par\q
\par We note that Corollary 4.1 is a case of [S6, Conjecture 4.3].
\par For example, let  $p$ be a prime such that $p\e 13\mod{24}$
and hence $p=c^2+d^2=x^2+3y^2$ with $c,d,x,y\in\Bbb Z$. Suppose $c\e
1\mod 4$, $d=2^rd_0$, $y=2^ty_0$ and $d_0\e y_0\e 1\mod 4$. If
$(c,x+d)=1$ or $(d_0,x+c)=1$, then
$$3^{\f{p-5}8}\e \cases \pm \f yx\mod p&\t{if $x\e \pm c\mod 3$,}
\\\mp\f{dy}{cx}\mod p&\t{if $x\e \pm d\mod 3$.}
\endcases$$
This partially solves [S5, Conjecture 9.1].

\pro{Theorem 4.2} Let $p$ and $q$ be primes such that $p\e 1\mod 4$,
$q\e 7\mod 8$, $p=c^2+d^2=x^2+qy^2$, $c,d,x,y\in\Bbb Z$, $c\e 1\mod
4$, $d=2^rd_0$, $y=2^ty_0$ and $d_0\e y_0\e 1\mod 4$.  Assume
$(c,x+d)=1$ or $(d_0,x+c)=1$. Suppose
 $\sls{c-di}{c+di}^{\f{q+1}8}\e i^m\mod q$.
 Then
$$q^{[p/8]}\e\cases (-1)^{\f y4}\sls dc^m\mod p&\t{if $p\e 1\mod
8$,}\\(-1)^{\f{x-1}2}\sls dc^m\f yx\mod p&\t{if $p\e 5\mod 8$.}
\endcases$$
\endpro
Proof. Observe that
$$\Ls{c-di}{c+di}^{\f{q+1}8}=\f{(c-di)^{\f{q+1}4}}{(c^2+d^2)^{\f{q+1}8}}
=\f{(c-di)^{\f{q+1}4}}{(x^2+qy^2)^{\f{q+1}8}}\e
\Ls{c-di}x^{\f{q+1}4}\mod q.$$ The result follows from Theorem 4.1.
\par\q
\par We note that if $q\nmid d$, then the $m$ in Theorem 4.2
depends only on $\f cd\mod q$.
 \pro{Corollary 4.2} Let $p\e 1\mod 4$ and $q\e 7\mod{8}$ be
primes such that $p=c^2+d^2=x^2+qy^2$ with $c,d,x,y\in\Bbb Z$ and
$q\mid cd(c^2-d^2)$. Suppose $c\e 1\mod 4$, $d=2^rd_0$, $y=2^ty_0$
and $d_0\e y_0\e 1\mod 4$. Assume $(c,x+d)=1$ or $(d_0,x+c)=1$.
\par $(\t{\rm i})$ If $p\e 1\mod 8$, then
$$q^{\f{p-1}8}\e \cases (-1)^{\f{q+1}8+\f y4}
\mod p&\t{if $q\mid c$,}
\\(-1)^{\f y4}\mod p&\t{if $q\mid d$,}
\\\pm (-1)^{\f{q+9}{16}+\f y4}\f dc\mod p&\t{if $16\mid (q-7)$ and
$c\e \pm d\mod q$,}
\\(-1)^{\f{q+1}{16}+\f y4}\mod p&\t{if $16\mid (q-15)$ and
$c\e \pm d\mod q$.}
\endcases$$
\par $(\t{\rm ii})$ If $p\e 5\mod 8$, then
$$q^{\f{p-5}8}\e \cases (-1)^{\f{q+1}8+\f{x-1}2} \f yx
\mod p&\t{if $q\mid c$,}\\(-1)^{\f{x-1}2}\f yx \mod p&\t{if $q\mid
d$,}\\\pm (-1)^{\f{q+9}{16}+\f{x-1}2}\f {dy}{cx}\mod p&\t{if $16\mid
(q-7)$ and $c\e \pm d\mod q$,}
\\(-1)^{\f{q+1}{16}+\f {x-1}2}\f yx\mod p&\t{if $16\mid (q-15)$ and
$c\e \pm d\mod q$.}
\endcases$$
\endpro
Proof. Clearly
$$\f{c-di}{c+di}\e \cases -1\mod q&\t{if $q\mid c$,}
\\1\mod q&\t{if $q\mid d$,}
\\-i\mod q&\t{if $c\e d\mod q$,}
\\i\mod q&\t{if $c\e -d\mod q$.}
\endcases$$
Thus the result follows from Theorem 4.2.
\par\q
\par In [S6] the author conjectured the condition $(c,x+d)=1$ or
$(d_0,x+c)=1$ in
Corollary 4.2 could be canceled.

\pro{Theorem 4.3} Let $p$ and $q$ be distinct primes of the form
$4k+1$, $p=c^2+d^2=x^2+qy^2$, $q=a^2+b^2$, $a,b,c,d,x,y\in\Bbb Z$,
$c\e 1\mod 4$, $d=2^rd_0$, $y=2^ty_0$ and $d_0\e y_0\e 1\mod 4$.
Assume $(c,x+d)=1$ or $(d_0,x+c)=1$. Suppose
$\sls{ac+bd}{ax}^{\f{q-1}4}\e \sls ba^m\mod q$.
\par $(\t{\rm i})$ If $p\e 1\mod 8$, then
$$q^{\f{p-1}8}
\e\cases (-1)^{\f d4+[\f {x+2}4]}\sls dc^m\mod p&\t{if $2\mid
x$,}\\(-1)^{\f{q-1}4\cdot\f{x-1}2+ \f d4+\f y4}\sls dc^m\mod p&\t{if
$2\nmid x.$}\endcases$$
\par $(\t{\rm ii})$ If $p\e 5\mod 8$, then
$$q^{\f{p-5}8}\e\cases
(-1)^{[\f x4]}\sls dc^{m+1}\f yx\mod p&\t{if $2\mid
x$,}\\-(-1)^{\f{q+3}4\cdot\f{x-1}2}\sls dc^{m+1}\f yx\mod p&\t{if
$2\nmid x$.}
\endcases$$
\endpro
Proof. Clearly $q\nmid x$. We first assume $(c,x+d)=1$. By Lemma
2.12 we have $(q,(x+d)(c^2+(x+d)^2))=1$. It is easily seen that
$\f{ac+b(x+d)}{ac-b(x+d)}\e \f{ac+bd}{ax}\cdot\f ba\mod q$.
 Thus, for $k=0,1,2,3$, using Lemma 2.7 we get
$$\align &\Qs{c+(x+d)i}q=i^k\\&\iff\f c{x+d}\in Q_k(q)
\iff \Ls{\f c{x+d}+\f ba}{\f c{x+d}-\f ba}^{\f{q-1}4}\e \Big(\f
ba\Big)^k\mod q
\\&\iff \Ls{ac+b(x+d)}{ac-b(x+d)}^{\f{q-1}4}\e \Ls ba^k\mod q
\\&\iff
\Big(\f{ac+bd}{ax}\cdot\f ba\Big)^{\f{q-1}4}\e \Ls ba^k\mod q
\\&\iff \Ls{ac+bd}{ax}^{\f{q-1}4}\e \Ls ba^{k-\f{q-1}4}\mod q.
\endalign$$
Since $ \sls{ac+bd}{ax}^{\f{q-1}4}\e \sls ba^m\mod q$, from the
above we get $\sqs{c+(x+d)i}q=i^{m+\f{q-1}4}$.  Now the result
follows from Theorems 3.1 and 3.2 immediately.
\par Suppose $(d_0,x+c)=1$. By Lemma 2.12 we
have $(q,(x+c)(d^2+(x+c)^2))=1$. It is easily seen that
$\f{ad-b(x+c)}{ad+b(x+c)}\e \f{ac+bd}{-ax}\mod q$.
 Thus, for $k=0,1,2,3$, using Lemma 2.7 we get
$$\align &\Qs{d-(x+c)i}q=i^k\\&\iff-\f d{x+c}\in Q_k(q)
\iff \Ls{-\f d{x+c}+\f ba}{-\f d{x+c}-\f ba}^{\f{q-1}4}\e \Big(\f
ba\Big)^k\mod q
\\&\iff \Ls{ad-b(x+c)}{ad+b(x+c)}^{\f{q-1}4}\e \Ls ba^k\mod q
\\&\iff
\Ls{ac+bd}{-ax}^{\f{q-1}4}\e \Ls ba^k\mod q
\\&\iff \Ls{ac+bd}{ax}^{\f{q-1}4}\e \Ls ba^{\f{q-1}2+k}\mod q.
\endalign$$
Since $ \sls{ac+bd}{ax}^{\f{q-1}4}\e \sls ba^m\mod q$, by the above
we get $\sqs{d-(x+c)i}q=i^{m-\f{q-1}2}$. Thus, applying Theorems
3.3-3.5 and the fact $\f x2\e \f{x^2}4=\f{c^2-qy^2+d^2}4 \e
\f{1-q}4+\f d2\mod 2$ in the case $2\mid x$  we derive the result.
The proof is now complete.

\pro{Corollary 4.3} Let $p\e 1\mod 4$ and $q\e 5\mod{8}$ be primes
such that $p=c^2+d^2=x^2+qy^2$ with $c,d,x,y\in\Bbb Z$ and $q\mid
cd$. Suppose $c\e 1\mod 4$, $d=2^rd_0$, $y=2^ty_0$ and $d_0\e y_0\e
1\mod 4$. Assume $(c,x+d)=1$ or $(d_0,x+c)=1$.
\par $(\t{\rm i})$ If $p\e 1\mod 8$, then
$$q^{\f{p-1}8}\e \cases \pm
(-1)^{\f d4+\f{x+2}4}\mod p&\t{if $2\mid x$ and $x\e \pm c\mod{q}$,}
\\\pm(-1)^{\f d4+\f{x-1}2+\f y4}\mod p&\t{if $2\nmid x$ and $x\e\pm
c\mod{q}$,}
\\\pm (-1)^{\f{q-5}8+\f d4+\f{x+2}4}\f dc\mod p
&\t{if $2\mid x$ and $x\e \pm d\mod{q}$,}
\\\pm (-1)^{\f{q-5}8+\f d4+\f{x-1}2+\f y4}\f dc\mod p
&\t{if $2\nmid x$ and $x\e \pm d\mod{q}$.}
\endcases$$
\par $(\t{\rm ii})$ If $p\e 5\mod 8$, then
$$q^{\f{p-5}8}\e \cases \pm \delta(x)\f {dy}{cx}\mod p&\t{if $x\e\pm
c\mod{q}$,}\\\mp (-1)^{\f{q-5 }8}\delta(x)\f yx\mod p&\t{if $x\e \pm
d\mod{q}$,}
\endcases$$
where $\delta(x)=1$ or $-1$ according as $8\mid x$ or not.
\endpro
Proof. If $x\e \pm c\mod q$, then $q\mid d$ and so
$\sls{ac+bd}{ax}^{\f{q-1}4}\e \sls cx^{\f{q-1}4}\e (\pm
1)^{\f{q-1}4} = \pm 1\mod q$. If $x\e \pm d\mod q$, then $q\mid c$
and so $\sls{ac+bd}{ax}^{\f{q-1}4}\e \sls {bd}{ax}^{\f{q-1}4}\e (\pm
\f ba)^{\f{q-1}4} \e \pm (-1)^{\f{q-5}8}\f ba\mod q$. Now putting
the above with Theorem 4.3 we deduce the result.

\pro{Theorem 4.4} Let $p$ and $q$ be distinct primes such that $p\e
1\mod 4$, $q\e 1\mod 8$, $p=c^2+d^2=x^2+qy^2$, $q=a^2+b^2$,
$a,b,c,d,x,y\in\Bbb Z$, $c\e 1\mod 4$, $d=2^rd_0$, $y=2^ty_0$ and
$d_0\e y_0\e 1\mod 4$. Assume $(c,x+d)=1$ or $(d_0,x+c)=1$. Suppose
$\sls{ac+bd}{ac-bd}^{\f{q-1}8}\e \sls ba^m\mod q$.
 \par $(\t{\rm i})$ If $p\e 1\mod 8$, then
$$q^{\f{p-1}8}\e (-1)^{\f d4+\f {xy}4}\Ls dc^m\mod p.$$

\par $(\t{\rm ii})$ If $p\e 5\mod 8$, then
$$q^{\f{p-5}8}\e\cases (-1)^{\f {x-2}4}\sls dc^{m+1}\f yx\mod p&\t{if $2\mid x$,}
\\(-1)^{\f{x+1}2}\sls dc^{m+1}\f yx\mod p&\t{if $2\nmid x$.}
\endcases$$
\endpro

Proof. Observe that $b^2\e -a^2\mod q$, $p\e x^2\mod q$ and so
$$\Ls{ac+bd}{ac-bd}^{\f{q-1}8}=\f{(ac+bd)^{\f{q-1}4}}{(a^2c^2-b^2d^2
)^{\f{q-1}8}} \e\f{(ac+bd)^{\f{q-1}4}}{(a^2p)^{\f{q-1}8}}\e
\Ls{ac+bd}{ax}^{\f{q-1}4}\mod q.$$ The result follows from Theorem
4.3.
\par\q
\par We note that if $q\nmid d$, then the $m$ in Theorem 4.4
depends only on $\f cd\mod q$. \pro{Corollary 4.4} Let $p\e 1\mod 4$
and $q\e 1\mod{8}$ be distinct primes such that $p=c^2+d^2=x^2+qy^2$
with $c,d,x,y\in\Bbb Z$ and $q\mid cd(c^2-d^2)$. Suppose $c\e 1\mod
4$, $d=2^rd_0$, $y=2^ty_0$ and $d_0\e y_0\e 1\mod 4$. Assume
$(c,x+d)=1$ or $(d_0,x+c)=1$.
\par $(\t{\rm i})$ If $p\e 1\mod 8$, then
$$q^{\f{p-1}8}\e \cases (-1)^{\f{q-1}8+\f d4+\f {xy}4}
\mod p&\t{if $q\mid c$,}
\\(-1)^{\f d4+\f {xy}4}\mod p&\t{if $q\mid d$,}
\\ (-1)^{\f{q-1}{16}+\f d4+\f {xy}4}\mod p
&\t{if $16\mid (q-1)$ and $c\e \pm d\mod q$,}
\\\pm (-1)^{\f{q-9}{16}+\f d4+\f {xy}4}\f dc\mod p
&\t{if $16\mid (q-9)$ and $c\e \pm d\mod q$.}
\endcases$$
\par $(\t{\rm ii})$ If $p\e 5\mod 8$, then
$$q^{\f{p-5}8}\e \cases (-1)^{\f{q-1}8+\f{x-2}4} \f {dy}{cx}
\mod p&\t{if $2\mid x$ and $q\mid c$,}
\\(-1)^{\f{q-1}8+\f{x+1}2} \f {dy}{cx}
\mod p&\t{if $2\nmid x$ and $q\mid c$,} \\(-1)^{\f{x-2}4} \f
{dy}{cx} \mod p&\t{if $2\mid x$ and $q\mid d$,}
\\(-1)^{\f{x+1}2}\f {dy}{cx} \mod
p&\t{if $2\nmid x$ and $q\mid d$,}\\ (-1)^{\f{q-1}{16}+\f{x-2}4}\f
{dy}{cx}\mod p&\t{if $2\mid x$, $16\mid (q-1)$ and $c\e \pm d\mod
q$,}
\\(-1)^{\f{q-1}{16}+\f{x+1}2}\f {dy}{cx}\mod p&\t{if $2\nmid x$,
$16\mid (q-1)$ and $c\e \pm d\mod q$,}
\\ \mp (-1)^{\f{q-9}{16}+\f{x-2}4}\f yx
\mod p&\t{if $2\mid x$, $16\mid (q-9)$ and $c\e \pm d\mod q$,}
\\\mp
(-1)^{\f{q-9}{16}+\f{x+1}2}\f yx\mod p&\t{if $2\nmid x$, $16\mid
(q-9)$ and $c\e \pm d\mod q$.}
\endcases$$
\endpro
Proof. Suppose that $q=a^2+b^2$ with $a,b\in\Bbb Z$. Then clearly
$$\Ls{ac+bd}{ac-bd}^{\f{q-1}8}\e \cases (-1)^{\f{q-1}8}\mod q&\t{if $q\mid c$,}
\\1\mod q&\t{if $q\mid d$,}
\\(-1)^{\f{q-1}{16}}\mod q&\t{if $16\mid (q-1)$ and $c\e \pm d\mod q$,}
\\\pm (-1)^{\f{q-9}{16}}\f ba\mod q&\t{if $16\mid (q-9)$ and $c\e \pm d\mod q$.}
\endcases$$
Thus the result follows from Theorem 4.4.

\pro{Corollary 4.5} Let $p\e 1\mod 4$ be a prime such that
$p\not=17$ and $p=c^2+d^2=x^2+17y^2$ with $c,d,x,y\in\Bbb Z$.
Suppose $c\e 1\mod 4$, $d=2^rd_0$, $y=2^ty_0$ and $d_0\e y_0\e 1\mod
4$. Assume $(c,x+d)=1$ or $(d_0,x+c)=1$.
\par $(\t{\rm i})$ If $p\e 1\mod 8$, then
$$17^{\f{p-1}8}\e \cases (-1)^{\f d4+\f {xy}4}
\mod p&\t{if $17\mid cd$,}
\\ -(-1)^{\f d4+\f {xy}4}\mod p
&\t{if  $c\e \pm d\mod {17}$,}
\\\pm (-1)^{\f d4+\f {xy}4}\f cd\mod p
&\t{if $c\e \pm 5d,\pm 10d\mod {17}$.}
\endcases$$
\par $(\t{\rm ii})$ If $p\e 5\mod 8$, then
$$17^{\f{p-5}8}\e \cases (-1)^{\f{x-2}4} \f {dy}{cx}
\mod p&\t{if $2\mid x$ and $17\mid cd$,}
\\(-1)^{\f{x+1}2} \f {dy}{cx}
\mod p&\t{if $2\nmid x$ and $17\mid cd$,} \\
(-1)^{\f{x+2}4}\f {dy}{cx}\mod p&\t{if $2\mid x$ and $c\e \pm d\mod
{17}$,}
\\(-1)^{\f{x-1}2}\f {dy}{cx}\mod p&\t{if $2\nmid x$ and $c\e \pm d\mod {17}$,}
\\ \pm (-1)^{\f{x-2}4}\f yx
\mod p&\t{if $2\mid x$ and $c\e \pm 5d,\pm 10d\mod {17}$,}
\\ \pm (-1)^{\f{x+1}2}\f yx
\mod p&\t{if $2\nmid x$ and $c\e \pm 5d,\pm 10d\mod {17}$.}
\endcases$$
\endpro
Proof. Since $17=1^2+4^2$ and $\ls{\pm 5+4}{\pm 5-4}^{\f{17-1}8} \e
\ls{\pm 10+4}{\pm 10-4}^{\f{17-1}8}\e \mp 4\mod{17}$, from  Theorem
4.4 and Corollary 4.4 we deduce the result.

\pro{Theorem 4.5} Let $p\e 1\mod 4$ be a prime,
$p=c^2+d^2=x^2+(a^2+b^2)y^2\not=a^2+b^2$, $a,b,c,d,x,y\in\Bbb Z$,
$a\not=0$, $2\mid a$, $(a,b)=1$, $c\e 1\mod 4$, $d=2^rd_0$,
$y=2^ty_0$ and $d_0\e y_0\e 1\mod 4$. Assume $(c,x+d)=1$ or
$(d_0,x+c)=1$. Suppose $\qs{(ac+bd)/x}{b+ai}=i^m$.
\par $(\t{\rm i})$ If $p\e 1\mod 8$, then
$$(a^2+b^2)^{\f{p-1}8}\e\cases (-1)^{\f d4+\f x4}\sls cd^m
\mod p&\t{if $4\mid a$ and $2\mid x$,}
\\(-1)^{\f d4+\f y4}\sls cd^m\mod p&\t{if $4\mid a$ and $2\nmid x$,}
\\(-1)^{\f{b+1}2+\f d4+\f {x-2}4}\sls cd^{m-1}
\mod p&\t{if $2\ \Vert\ a$ and $2\mid x$,}
\\(-1)^{\f{b-1}2+\f d4+\f y4+\f{x-1}2}\sls cd^{m-1}\mod p
&\t{if $2\ \Vert\ a$ and $2\nmid x$.}
\endcases$$
\par $(\t{\rm ii})$ If $p\e 5\mod 8$, then
$$(a^2+b^2)^{\f{p-5}8}\e\cases
(-1)^{\f {x-2}4}\sls cd^{m-1}\f yx\mod p&\t{if $4\mid a$ and $2\mid
x$,}\\(-1)^{\f{x+1}2}\sls cd^{m-1}\f yx\mod p&\t{if $4\mid a$ and
$2\nmid x$,}
\\(-1)^{\f x4+\f{b+1}2}\sls cd^m\f yx\mod p&\t{if $2\ \Vert\ a$ and $2\mid
x$,}\\(-1)^{\f{b-1}2}\sls cd^m\f yx\mod p&\t{if $2\ \Vert\ a$ and
$2\nmid x$.}
\endcases$$
\endpro
Proof. Set $q=a^2+b^2$.  Then clearly $2\nmid q$ and $p\nmid q$.
 We first assume $(c,x+d)=1$.
By Lemma 2.12 we have $(q,x+d)=(q,c^2+(x+d)^2)=1$. Since
$\f{c-(x+d)i}{c+(x+d)i}\e \f{c-di}{ix}\mod q$,  we see that
$$\aligned &\Qs{c/(x+d)+i}q\\& =\Qs{c+(x+d)i}q
=\Qs{c+(x+d)i}{b+ai}\Qs{c+(x+d)i}{b-ai}
\\&=\Qs{c+(x+d)i}{b+ai}\overline{\Qs{c-(x+d)i}{b+ai}}
=\Qs{c+(x+d)i}{b+ai}\Qs{c-(x+d)i}{b+ai}^{-1}
\\&=\Qs{\f{c-(x+d)i}{c+(x+d)i}}{b+ai}^{-1}
=\Qs{\f{c-di}{ix}}{b+ai}^{-1}=\Qs{ai}{b+ai}
\Qs{(ac-adi)/x}{b+ai}^{-1}
\\&=\Qs {-b}{b+ai}\Qs{(ac+bd)/x}
{b+ai}^{-1}=(-1)^{\f {b+1}2\cdot \f a2}\Qs {b+ai}bi^{-m}
\\&=(-1)^{\f {b+1}2\cdot \f a2}\Qs ibi^{-m}
=(-1)^{\f {b+1}2\cdot \f a2+\f{b^2-1}8}i^{-m} \\&= (-1)^{\f
{b+1}2\cdot \f a2+\f{q-1-a^2}8}i^{-m}=(-1)^{\f {b+1}2\cdot \f a2+[\f
q8]}i^{-m}.
\endaligned$$
 This together with
Theorems 3.1 and 3.2 yields the result under the condition
$(c,x+d)=1$.

\par Now we assume $(d_0,x+c)=1$. By Lemma 2.12 we have
$(q,x+c)=(q,(x+c)^2+d^2)=1$. Since $\f{d+(x+c)i}{d-(x+c)i}\e
\f{c-di}{-x}\mod q$, using Lemma 2.6 we see that
$$\align &\Qs{d/(x+c)-i}q\\& =\Qs{d-(x+c)i}q
=\Qs{d-(x+c)i}{b+ai}\Qs{d-(x+c)i}{b-ai}
\\&=\Qs{d-(x+c)i}{b+ai}\overline{\Qs{d+(x+c)i}{b+ai}}
=\Qs{d-(x+c)i}{b+ai}\Qs{d+(x+c)i}{b+ai}^{-1}
\\&=\Qs{\f{d+(x+c)i}{d-(x+c)i}}{b+ai}^{-1}
=\Qs{\f{c-di}{-x}}{b+ai}^{-1}=\Qs{-a}{b+ai}
\Qs{(ac-adi)/x}{b+ai}^{-1}
\\&=(-1)^{\f a2}\Qs{a}{b+ai}\Qs{(ac+bd)/x}{b+ai}^{-1}
=(-1)^{\f a2}\Qs{a}{b+ai}i^{-m} \\&=\cases -(-1)^{\f{b-1}2}i\cdot
i^{-m}=(-1)^{\f{b+1}2}i^{1-m}&\t{if $2\ \Vert\ a$,}\\1 \cdot 1\cdot
i^{-m}=i^{-m}&\t{if $4\mid a$.}
\endcases
\endalign$$
 Combining this  with
Theorems 3.3-3.5 we deduce the result under the condition
$(d_0,x+c)=1$. The proof is now complete.
\par\q
\newline{\bf Remark 4.1}  We conjecture that the condition $(c,x+d)=1$ or
$(d_0,x+c)=1$ in Theorems 4.1-4.5 and Corollaries 4.1-4.5 can be
canceled. See [S5] and [S6].
 \subheading{5.
Congruences for $\ls{b+\sqrt{b^2+2^{2\a}}}2^{\f{p-1}4}\mod p$}
\pro{Theorem 5.1} Let $p\e 1\mod 4$ be a prime, $b\in\Bbb Z$,
$2\nmid b$, $p\not=b^2+4$ and $p=c^2+d^2=x^2+(b^2+4)y^2$ for some
$c,d,x,y\in\Bbb Z$. Suppose $c\e 1\mod 4$, $d=2^rd_0$,
$x=2^sx_0(2\nmid x_0)$, $y=2^ty_0$ and $d_0\e y_0\e 1\mod 4$. Assume
$(c,x+d)=1$ or $(d_0,x+c)=1$.
\par $(\t{\rm i})$ If $4\nmid xy$, then
$$\aligned\Ls{b+(-1)^{\f{x_0-1}2}\f{cx}{dy}}2^{\f{p-1}4}&\e
-\Ls{b-(-1)^{\f{x_0-1}2}\f{cx}{dy}}2^{\f{p-1}4}
\\&\e\cases (-1)^{[\f b4]+\f d4}\f cd\mod p&\t{if $2\ \Vert\ x$,}
\\ 1\mod p&\t{if $2\ \Vert\ y$.}
\endcases\endaligned$$
\par $(\t{\rm ii})$ If $4\mid xy$, then
$$\aligned\Ls{b+\f{cx}{dy}}2^{\f{p-1}4}\e \Ls{b-\f{cx}{dy}}2^{\f{p-1}4}
\e \cases (-1)^{[\f b4]+\f x4}\f cd\mod p&\t{if $4\mid x$,}
\\(-1)^{\f d4+\f y4} \mod p&\t{if $4\mid y$.}
\endcases\endaligned$$
\endpro
Proof. Since $\ls{-1}{b+2i}_4=-1$, by [S5, Corollary 6.1] we have
$$\aligned&\Ls{b-(-1)^{\f{x_0-1}2}cx/(dy)}2^{\f{p-1}4}\\&\e\cases\mp
(-1)^{\f{b-1}2} (b^2+4)^{\f{p-5}8}\f xy \mod p &\t{if $2\ \Vert\ y$
and $\qs {(2c+bd)/x}{b+2i}=\pm 1$,}
\\\mp
(-1)^{\f{b-1}2} (b^2+4)^{\f{p-5}8}\f {dx}{cy} \mod p &\t{if $2\
\Vert\ y$ and $\qs {(2c+bd)/x}{b+2i}=\pm i$,}
\\\mp (-1)^{\f{b-1}2+\f{x_0-1}2}(b^2+4)^{\f{p-1}8}\f dc\mod p&\t{if $4\mid y$
and $\qs {(2c+bd)/x}{b+2i}=\pm 1$,}
\\\pm (-1)^{\f{b-1}2+\f{x_0-1}2}(b^2+4)^{\f{p-1}8}\mod p&\t{if $4\mid y$
and $\qs {(2c+bd)/x}{b+2i}=\pm i$,}
\\\mp
(-1)^{\f{b^2-1}8+\f{x_0-1}2} (b^2+4)^{\f{p-1}8} \mod p &\t{if $2\
\Vert\ x$ and $\qs {(2c+bd)/x}{b+2i}=\pm 1$,}
\\\mp (-1)^{\f{b^2-1}8+\f{x_0-1}2}(b^2+4)^{\f{p-1}8}\f dc\mod p&\t{if $2\ \Vert\ x$
and $\qs {(2c+bd)/x}{b+2i}=\pm i$,}
\\\pm
(-1)^{\f{b^2-1}8} (b^2+4)^{\f{p-5}8}\f {dx}{cy} \mod p &\t{if $4\mid
x$ and $\qs {(2c+bd)/x}{b+2i}=\pm 1$,}
\\\mp (-1)^{\f{b^2-1}8}(b^2+4)^{\f{p-5}8}\f xy\mod p&\t{if $4\mid x$
and $\qs {(2c+bd)/x}{b+2i}=\pm i$.}
\endcases\endaligned\tag 5.1$$
Taking $a=2$ in Theorem 4.5 we obtain
$$(b^2+4)^{[\f p8]}\e\cases (-1)^{\f{b+1}2+\f d4+\f{x_0-1}2}\sls
cd^{m-1} \mod p&\t{if $2\ \Vert\ x$,}
\\(-1)^{\f{b-1}2}\sls cd^m\f yx\mod p&\t{if $2\ \Vert\ y$,}
\\(-1)^{\f x4+\f{b+1}2}\sls cd^m\f yx\mod p&\t{if $4\mid x$,}
\\(-1)^{\f{b-1}2+\f d4+\f y4+\f{x_0-1}2}\sls cd^{m-1}\mod p&\t{if
$4\mid y$.}\endcases$$
This together with (5.1) yields
$$\Ls{b-(-1)^{\f{x_0-1}2}\f{cx}{dy}}2^{\f{p-1}4}
\e\cases -(-1)^{[\f b4]+\f d4}\f cd\mod p&\t{if $2\ \Vert\ x$,}
\\-1\mod p&\t{if $2\ \Vert\ y$,}
\\(-1)^{[\f b4]+\f x4}\f cd\mod p&\t{if $4\mid x$,}
\\(-1)^{\f d4+\f y4}\mod p&\t{if $4\mid y$.}
\endcases\tag 5.2$$
Since $\f{b+(-1)^{\f{x_0-1}2}\f{cx}{dy}}2\cdot \f
{b-(-1)^{\f{x_0-1}2}\f{cx}{dy}}2\e \f{b^2-(b^2+4)}4=-1\mod p$, we
see that
$$\Ls{b\pm (-1)^{\f{x_0-1}2}\f{cx}{dy}}2^{\f{p-1}4}\e
(-1)^{\f{p-1}4}\Ls{b\mp (-1)^{\f{x_0-1}2}\f{cx}{dy}}2^{-\f{p-1}4}
\mod p.\tag 5.3$$ Now combining (5.2) with (5.3) we deduce the
result.

 \pro{Corollary 5.1}
Let $p\e 1\mod 4$ be a prime, $b\in\Bbb Z$, $2\nmid b$,
$p\not=b^2+4$ and $p=c^2+d^2=x^2+(b^2+4)y^2$ with $c,d,x,y\in\Bbb Z$
and $4\mid xy$. Suppose $c\e 1\mod 4$, $d=2^rd_0$ and $d_0\e 1\mod
4$. Assume $(c,x+d)=1$ or $(d_0,x+c)=1$. Then
$$\Ls{b+\sqrt{b^2+4}}2^{\f{p-1}4}
\e \cases (-1)^{[\f b4]+\f x4}\f cd\mod p&\t{if $4\mid x$,}
\\(-1)^{\f d4+\f y4} \mod p&\t{if $4\mid y$.}
\endcases$$\endpro
\par\q
\newline{\bf Remark 5.1} Let $p$ be a prime of the form $4k+1$
and let $b$ be an odd number such that $\sls{b^2+4}p=1$. For the
congruences for $\sls{b+\sqrt{b^2+4}}2^{\f{p-1}2}\mod p$, one may
consult [Lem], [S3] and [S4].

 \pro{Theorem 5.2}
Let $p\e 1\mod 4$ be a prime, $\a\in\{2,3,4,\ldots\}$, $b\in\Bbb Z$,
 $2\nmid b$, $p\not=2^{2\a}+b^2$ and
$p=c^2+d^2=x^2+(2^{2\a}+b^2)y^2$ for some $c,d,x,y\in\Bbb Z$.
Suppose $c\e 1\mod 4$, $d=2^rd_0$, $x=2^sx_0(2\nmid x_0)$,
$y=2^ty_0$ and $d_0\e y_0\e 1\mod 4$. Assume $(c,x+d)=1$ or
$(d_0,x+c)=1$.
\par $(\t{\rm i})$ If $p\e 1\mod 8$, then
$$\aligned\Ls{b-\f{cx}{dy}}2^{\f{p-1}4}\e
\Ls{b+\f{cx}{dy}}2^{\f{p-1}4}\e\cases
(-1)^{\f{b^2-1}8+2^{\a-2}+\f{d+x}4\a}\mod p&\t{if $4\mid x$,}
\\(-1)^{\f {d+y}4\a}\mod p&\t{if $4\mid y$.}
\endcases\endaligned$$
\par $(\t{\rm ii})$ If $p\e 5\mod 8$, then
$$\aligned\Ls{b-(-1)^{\f{x_0-1}2}\f{cx}{dy}}2^{\f{p-1}4}&\e (-1)^{\a}
\Ls{b+(-1)^{\f{x_0-1}2}\f{cx}{dy}}2^{\f{p-1}4}
\\&\e \cases (-1)^{\f{(b+2)^2-1}8+\f{b+1}2\a}
\mod p&\t{if $2\ \Vert\ x$,}
\\(-1)^{\f{b-1}2(\a+1)}\mod p&\t{if $2\ \Vert\ y$.}
\endcases\endaligned$$
\endpro
Proof.  It is clear that
$$\aligned&\Ls{b+\f{cx}{dy}}2^{\f{p-1}4}\Ls{b-\f{cx}{dy}}2^{\f{p-1}4}
\\&\e\Ls{b^2-(2^{2\a}+b^2)}4^{\f{p-1}4}=(-1)^{\f{p-1}4}2^{\f{p-1}2(\a-1)}
\e (-1)^{\f{p-1}4\a}\mod p.\endaligned\tag 5.4$$  By [IR, p.64] we
have
$$2^{\f{p-1}4}\e \Ls dc^{\f{cd}2}\e \Ls dc^{\f d2}\e
\cases (-1)^{\f d4}\mod p&\t{if $p\e 1\mod 8$,}
\\\f dc\mod p&\t{if $p\e 5\mod 8$.}\endcases\tag 5.5$$
\par Suppose that $\sls{(2^{\a}c+bd)/x}{b+2^{\a}i}_4=i^m$.
As $\sls{-1}{b+2^{\a}i}_4=1$ we have
$\sls{(2^{\a}c+bd)/(-x)}{b+2^{\a}i}_4=i^m$. If $p\e 1\mod 8$ and
$4\mid x$, by [S5, Theorem 6.1], (5.5)
  and Theorem 4.5 we have
$$\aligned\Ls{b-(-1)^{\f{x_0-1}2}\f{cx}{dy}}2^{\f{p-1}4}
&\e (-1)^{\f{b^2-1}8+(\a+1)\f{2^{\a}-(-1)^{\f{x_0-1}2}x}4+(\a-1)\f
d4}\Ls dc^m(2^{2\a}+b^2)^{\f{p-1}8}
\\&\e (-1)^{\f{b^2-1}8+(\a+1)\f{2^{\a}-(-1)^{\f{x_0-1}2}x}4+(\a-1)\f d4}\Ls
dc^m\cdot (-1)^{\f d4+\f x4}\Ls cd^m
\\&=(-1)^{\f{b^2-1}8+2^{\a-2}(\a+1)+\f d4\a+\f x4\a}
=(-1)^{\f{b^2-1}8+2^{\a-2}+\f {d+x}4\a}\mod p.\endaligned$$ If $p\e
1\mod 8$ and $4\mid y$, by [S5, Theorem 6.1], (5.5)
  and Theorem 4.5 we have
$$\aligned\Ls{b-(-1)^{\f{x_0-1}2}\f{cx}{dy}}2^{\f{p-1}4}
&\e (-1)^{(\a+1)\f y4+(\a-1)\f d4}\Ls dc^m(2^{2\a}+b^2)^{\f{p-1}8}
\\&\e (-1)^{(\a+1)\f y4+(\a-1)\f d4}\Ls dc^m
\cdot (-1)^{\f d4+\f y4}\Ls cd^m
\\&=(-1)^{\f {d+y}4\a}\mod p.\endaligned$$
If $p\e 5\mod 8$ and $2\ \Vert\ y$, by [S5, Theorem 6.1], (5.5)
  and Theorem 4.5 we have
$$\aligned&\Ls{b-(-1)^{\f{x_0-1}2}\f{cx}{dy}}2^{\f{p-1}4}
\\&\e (-1)^{\f{b+1}2}\Ls dc^{m-b\a+\a-1}(-1)^{\f{x_0-1}2}\f
xy(2^{2\a}+b^2)^{\f{p-5}8}
\\&\e (-1)^{\f{b+1}2}\Ls dc^{m-b\a+\a-1}(-1)^{\f{x_0-1}2}\f xy
\cdot (-1)^{\f{x_0-1}2}\Ls cd^{m+1}\f yx
\\&\e (-1)^{\f{b-1}2(\a+1)}
\mod p.\endaligned$$ If $p\e 5\mod 8$ and $2\ \Vert\ x$, by [S5,
Theorem 6.1], (5.5)
  and Theorem 4.5 we have
$$\aligned&\Ls{b-(-1)^{\f{x_0-1}2}\f{cx}{dy}}2^{\f{p-1}4}
\\&\e -(-1)^{\f{(b+2)^2-9}8}\Ls
dc^{m+(-1)^{\f{b-1}2}\a+\a-1}(-1)^{\f{x_0-1}2}\f
xy(2^{2\a}+b^2)^{\f{p-5}8}
\\&\e -(-1)^{\f{(b+2)^2-9}8}\Ls dc^{m+(-1)^{\f{b-1}2}\a+\a-1}(-1)^{\f{x_0-1}2}
\f xy \cdot (-1)^{\f{x-2}4}\Ls cd^{m-1}\f yx
\\&\e (-1)^{\f{(b+2)^2-1}8+\f{b+1}2\a}
\mod p.\endaligned$$ Now putting all the above together  we derive
the result.

\par\q
\par We remark that $(-1)^{\f{(a_0+2)^2-(b+2)^2}8}$ should be
$-(-1)^{\f{(a_0+2)^2-(b+2)^2}8}$ in [S5, Theorem 6.1(ii)].
\par Taking $\a=2$ in Theorem 5.2 we deduce the following result.
 \pro{Corollary 5.2} Let $p\e 1\mod 4$ be a prime,
$b\in\Bbb Z$, $2\nmid b$, $p\not=b^2+16$ and
$p=c^2+d^2=x^2+(b^2+16)y^2$ for some $c,d,x,y\in\Bbb Z$. Suppose
$c\e 1\mod 4$ and $d=2^rd_0$ with $d_0\e 1\mod 4$. Assume
$(c,x+d)=1$ or $(d_0,x+c)=1$. Then
$$\Ls{b+\sqrt{b^2+16}}2^{\f{p-1}4}\e
\cases (-1)^y\mod p&\t{if $b\e 1\mod 8$,}
\\(-1)^{\f{p-1}4}\mod p&\t{if $b\e 3\mod 8$,}
\\1\mod p&\t{if $b\e 5\mod 8$,}
\\(-1)^{\f{p-1}4+y}\mod p&\t{if $b\e 7\mod 8$.}
\endcases$$
\endpro
\par\q
\newline
{\bf Remark 5.2} Whether $(c,x+d)=1$ or not, the author conjectured
the congruences in Theorem 5.1 and so Corollary 5.1. See [S5,
Conjecture 9.5]. We conjecture that the condition $(c,x+d)=1$ or
$(d_0,x+c)=1$ in Theorem 5.2 and Corollary 5.2 can also be canceled.

\subheading{6. Congruences for $U_{\f{p-1}4}(b,-4^{\a-1})$ and
$V_{\f{p-1}4}(b,-4^{\a-1})\mod p$}

\par Let $\{U_n(b,c)\}$ and $\{V_n(b,c)\}$ be the Lucas sequences given in Section 1.
Set $d=b^2-4c$. It is well known that for $n\in\Bbb N$,
$$U_n(b,c)=\cases \f 1{\sqrt d}\big\{\ls{b+\sqrt d}2^n-
\ls{b-\sqrt d}2^n\big\}&\t{if $d\not=0$,}
\\n(\f b2)^{n-1}&\t{if $d=0$}\endcases\tag 6.1$$
and
$$V_n(b,c)=\Ls {b+\sqrt d}2^n+
\Ls{b-\sqrt d}2^n.\tag 6.2$$ From [S1, Lemma 6.1(b)] we know that if
$p>3$ is a prime such that $p\nmid bcd$, then
$$p\mid U_n(b,c)\iff V_{2n}(b,c)\e
2c^n\mod p.\tag 6.3$$

\pro{Theorem 6.1} Let $p\e 1\mod 4$ be a prime, $b\in\Bbb Z$,
$2\nmid b$ and $p= c^2+d^2=x^2+(b^2+4)y^2\not=b^2+4$ for some
$c,d,x,y\in\Bbb Z$. Suppose $c\e 1\mod 4$, $d=2^rd_0$, $y=2^ty_0$
and $d_0\e y_0\e 1\mod 4$. Assume $(c,x+d)=1$ or $(d_0,x+c)=1$.
\par $(\t{\rm i})$ If $4\nmid xy$, then $V_{\f{p-1}4}(b,-1)\e 0\mod
p$ and
$$U_{\f{p-1}4}(b,-1)\e \cases(-1)^{[\f b4]+\f d4+\f{x-2}4}\f{2y}x
\mod p&\t{if $2\ \Vert\ x$,}
\\(-1)^{\f{x-1}2}\f{2dy}{cx}\mod p&\t{if $2\ \Vert\ y$.}
\endcases$$

\par $(\t{\rm ii})$ If $4\mid xy$, then
 $U_{\f{p-1}4}(b,-1)\e 0\mod
p$ and
$$V_{\f{p-1}4}(b,-1)\e \cases
2(-1)^{[\f b4]+\f x4}\f cd\mod p&\t{if $4\mid x$,}
\\2(-1)^{\f {d+y}4}
\mod p&\t{if $4\mid y$.}
\endcases$$
\endpro
Proof. Suppose $x=2^sx_0$ with $2\nmid x_0$. Since $\sls{cx}{dy}^2\e
b^2+4\mod p$, using (6.1) and (6.2) we see that
$$U_{\f{p-1}4}(b,-1)\e (-1)^{\f{x_0-1}2} \f{dy}{cx}
  \Big\{\Ls{b+(-1)^{\f{x_0-1}2}\f{cx}{dy}}2^{\f{p-1}4}
-\Ls{b-(-1)^{\f{x_0-1}2}\f{cx}{dy}}2^{\f{p-1}4}\Big\}\mod p$$ and
$$V_{\f{p-1}4}(b,-1)\e \Ls{b+\f{cx}{dy}}2^{\f{p-1}4}
+\Ls{b-\f{cx}{dy}}2^{\f{p-1}4}\mod p.$$ Now applying Theorem 5.1 we
deduce the result.
\par\q
\newline{\bf Remark 6.1} When $p$ is a prime of the form $8k+1$ and
$p=c^2+d^2=x^2+(b^2+4)y^2$ with $b\in\{1,3\}$ and $4\mid y$, the
conjecture $V_{\f{p-1}4}(b,-1)\e 2(-1)^{\f{d+y}4}\mod p$ is
equivalent to a conjecture of E. Lehmer. See [L, Conjecture 4].
 \pro{Theorem 6.2} Let $p\e 1\mod 8$ be a prime,
$b\in\Bbb Z$, $2\nmid b$ and $p=c^2+d^2=x^2+(b^2+4)y^2\not=b^2+4$
for some $c,d,x,y\in\Bbb Z$. Suppose $c\e 1\mod 4$ and $d=2^rd_0$
with $d_0\e 1\mod 4$. Assume $(c,x+d)=1$ or $(d_0,x+c)=1$. Then
$p\mid U_{\f{p-1}8}(b,-1)$ if and only if $y\e \f{p-1}2+d\mod 8$.
\endpro
Proof. This is immediate from (6.3) and Theorem 6.1.
\par\q
\par Using (6.1), (6.2) and Theorem 5.2 one can similarly prove the
following result.
 \pro{Theorem 6.3}
Let $p\e 1\mod 4$ be a prime, $b,\a\in\Bbb Z$, $\a\ge 2$, $2\nmid
b$, $p\not=b^2+4^{\a}$ and $p=c^2+d^2=x^2+(b^2+4^{\a})y^2$ for some
$c,d,x,y\in\Bbb Z$. Suppose $c\e 1\mod 4$, $d=2^rd_0$, $y=2^ty_0$
and $d_0\e y_0\e 1\mod 4$. Assume $(c,x+d)=1$ or $(d_0,x+c)=1$.
\par $(\t{\rm i})$ If $p\e 1\mod 8$, then
$U_{\f{p-1}4}(b,-4^{\a-1})\e 0\mod p$ and
$$V_{\f{p-1}4}(b,-4^{\a-1})\e
 \cases 2(-1)^{\f{b^2-1}8+2^{\a-2}+\f{d+x}4\a}\mod p&\t{if $4\mid
x$,} \\2(-1)^{\f {d+y}4\a}\mod p&\t{if $4\mid y$.}
\endcases$$
\par $(\t{\rm ii})$ If $p\e 5\mod 8$, then
$$U_{\f{p-1}4}(b,-4^{\a-1})\e \cases 0\mod p&\t{if $2\mid \a$,}
\\2(-1)^{\f{(b+2)^2-1}8+\f{b-1}2+\f{x-2}4}\f{dy}{cx}\mod p&\t{if
$2\nmid \a$ and $2\ \Vert\ x$,}
\\2(-1)^{\f{x+1}2}\f{dy}{cx}\mod p&\t{if
$2\nmid \a$ and $2\ \Vert\ y$}
\endcases$$
and
$$V_{\f{p-1}4}(b,-4^{\a-1})\e \cases 0\mod p&\t{if $2\nmid \a$,}
\\2(-1)^{\f{(b+2)^2-1}8}\mod p&\t{if
$2\mid \a$ and $2\ \Vert\ x$,}
\\2(-1)^{\f{b-1}2}\mod p&\t{if
$2\mid \a$ and $2\ \Vert\ y$.}
\endcases$$
\endpro

\par From (5.5), (6.3) and Theorem 6.3 we derive the following result.
\pro{Theorem 6.4} Let $p\e 1\mod 8$ be a prime, $b,\a\in\Bbb Z$,
$\a\ge 2$, $2\nmid b$, $p\not=b^2+4^{\a}$ and
$p=c^2+d^2=x^2+(b^2+4^{\a})y^2$ for some $c,d,x,y\in\Bbb Z$. Suppose
$c\e 1\mod 4$ and $d=2^rd_0$ with $d_0\e 1\mod 4$. Assume
$(c,x+d)=1$ or $(d_0,x+c)=1$. Then
$$p\mid
U_{\f{p-1}8}(b,-4^{\a-1})\iff \f{p-1}8+\f d4\e
 \cases \f{b^2-1}8+2^{\a-2}+\f{x}4\a\mod 2&\t{if $4\mid
x$,} \\\f {y}4\a\mod 2&\t{if $4\mid y$.}
\endcases$$
\endpro
\subheading{7. Congruences for $U_{\f{p-1}4}(4a,-1)$ and
$V_{\f{p-1}4}(4a,-1)\mod p$}

\pro{Theorem 7.1} Let $a\in\Bbb Z$, and let $p\e 1\mod 4$ be a prime
such that $p=c^2+d^2=x^2+(4a^2+1)y^2$ with $c,d,x,y\in\Bbb Z$, $c\e
1\mod 4$ and $p\not=4a^2+1$. Suppose $d=2^rd_0$, $y=2^{\beta}y_0$
and $d_0\e y_0\e 1\mod 4$. Assume that $(c,x+d)=1$ or $(d_0,x+c)=1$.
\par $(\t{\rm i})$ If $p\e 1\mod 8$, then
$$U_{\f{p-1}4}(4a,-1)\e\cases (-1)^{\f{a+1}2+\f d4+\f{x-2}4}\f yx\mod p&\t{if $2\nmid ay$,}
\\0\mod p&\t{if $2\mid ay$}\endcases $$
and
$$V_{\f{p-1}4}(4a,-1)\e\cases 2(-1)^{\f d4+\f a2y+\f{xy}4}\mod
p&\t{if $2\mid a$,}
\\2(-1)^{\f d4+\f y4}\mod p&\t{if $2\nmid a$ and $2\mid y$,}
\\0\mod p&\t{if $2\nmid ay$.}\endcases$$

\par $(\t{\rm ii})$ If $p\e 5\mod 8$, then
$$U_{\f{p-1}4}(4a,-1)\e\cases (-1)^{\f a2+\f{x-2}4}\f {dy}{cx}\mod p
&\t{if $2\mid a$ and $2\nmid y$,}
\\(-1)^{\f{x+1}2}\f{dy}{cx}\mod p&\t{if $2\mid a$ and $2\mid y$,}
\\\f{dy}{cx}\mod p&\t{if $2\nmid a$ and $2\mid y$,}
\\0\mod p&\t{if $2\nmid ay$}\endcases $$
and
$$V_{\f{p-1}4}(4a,-1)\e\cases 0\mod
p&\t{if $2\mid ay$,}
\\2(-1)^{\f {a-1}2+\f x4}\f dc\mod p&\t{if $2\nmid ay$.}\endcases$$

\endpro
Proof. We first assume $p\e 1\mod 8$. Since
$\sqs{-1}{1-2ai}=(-1)^a$, by [S5, Theorem 7.3] we have
$$\aligned &U_{\f{p-1}4}(4a,-1)\\&\e\cases
(-1)^{\f{a-1}2}\f{dy}{cx}(-\f dc)^k(4a^2+1)^{\f{p-1}8}\mod p&\t{if
$2\nmid ay$ and $\sqs{(-d-2ac)/x}{1-2ai}=i^k$,}
\\0\mod p&\t{if $2\mid ay$}
\endcases\endaligned$$ and
$$\aligned&V_{\f{p-1}4}(4a,-1)\\&\e \cases
 2(-1)^{\f a2y}(-\f dc)^k(4a^2+1)^{\f{p-1}8}\mod p&\t{if $2\mid a$ and
$\qs{(-d-2ac)/x}{1-2ai}=i^k$,}
\\ 2(-1)^{\f{x-1}2}\f dc(-\f dc)^k(4a^2+1)^{\f{p-1}8}\mod p&\t{if
$2\nmid a$, $2\mid y$ and $\sqs{(-d-2ac)/x}{1-2ai}=i^k$,}
\\0\mod p&\t{if $2\nmid ay.$}
\endcases\endaligned$$
If $(c,x+d)=1$, then $(4a^2+1,(x+d)(c^2+(x+d)^2))=1$ by Lemma 2.12.
It is easily seen that $\f{c+(x+d)i}{c-(x+d)i}\e
\f{-d+ic}x\mod{4a^2+1}.$ Thus,
$$\align
&\Qs{c/(x+d)+i}{4a^2+1}\\&=\Qs{c+(x+d)i}{4a^2+1}=\Qs{c+(x+d)i}{1-2ai}\Qs{c+(x+d)i}{1+2ai}
\\&=\Qs{c+(x+d)i}{1-2ai}\Qs{c-(x+d)i}{1-2ai}^{-1}=\Qs{(c+(x+d)i)/(c-(x+d)i)}{1-2ai}
\\&=\Qs{(-d+ic)/x}{1-2ai}=\Qs{(-d-2ac)/x}{1-2ai}.\endalign$$
If $(d_0,x+c)=1$,  as $(4a^2+1)y^2=(x+c)^2-2x(x+c)+d^2$ we have
$(x+c,(4a^2+1)y^2)=(x+c,d^2)$ and so $(4a^2+1,x+c)=1$. Also,
$(4a^2+1,d^2+(x+c)^2)=(4a^2+1,
(4a^2+1)y^2+2x(x+c))=(4a^2+1,2x(x+c))=1$. It is easily seen that
$\f{d-(x+c)i}{d+(x+c)i}\e -\f{c+di}x\mod{4a^2+1}$. Thus,

$$\align
&\Qs{d/(x+c)-i}{4a^2+1}\\&=\Qs{d-(x+c)i}{4a^2+1}=\Qs{d-(x+c)i}{1-2ai}\Qs{d-(x+c)i}{1+2ai}
\\&=\Qs{d-(x+c)i}{1-2ai}\Qs{d+(x+c)i}{1-2ai}^{-1}=\Qs{(d-(x+c)i)/(d+(x+c)i)}{1-2ai}
\\&=\Qs{-(c+di)/x}{1-2ai}=\Qs i{1-2ai}\Qs{(-d+ci)/x}{1-2ai}
\\&=i^{(1-(-1)^a)/2}\Qs{(-d-2ac)/x}{1-2ai}.\endalign$$
By Theorems 3.1 and 3.2, if $(c,x+d)=1$ and
$\sqs{c/(x+d)+i}{4a^2+1}=i^k$, then
$$(4a^2+1)^{\f{p-1}8}
\e\cases (-1)^{\f d4+\f {xy}4}\sls dc^k\mod p&\t{if $2\mid a$,}
\\(-1)^{\f d4+\f{x-2}4}\sls dc^{k+1}\mod p&\t{if $2\nmid ay$,}
\\(-1)^{\f d4+\f{x-1}2+\f y4}\sls dc^{k-1}\mod p&\t{if $2\nmid a$ and $2\mid y$.}
\endcases$$
By Theorems 3.3-3.5, if $(d_0,x+c)=1$ and
$\sqs{d/(x+c)-i}{4a^2+1}=i^k$, then
$$(4a^2+1)^{\f{p-1}8}\e\cases (-1)^{\f d4+[\f x4]}\sls dc^k\mod
p&\t{if $2\mid x$,}
\\(-1)^{a\cdot\f{x+1}2+\f d4+\f y4}\sls dc^k\mod p&\t{if $2\nmid
x$.}\endcases$$
 Now combining all the above we deduce (i). Part (ii) can be proved
 similarly.
 \pro{Corollary 7.1} Let $a\in\Bbb Z$, and let $p\e 1\mod 4$ be a prime
such that $p=c^2+d^2=x^2+(4a^2+1)y^2$ with $c,d,x,y\in\Bbb Z$, $c\e
1\mod 4$ and $p\not=4a^2+1$. Suppose $d=2^rd_0$ with $d_0\e 1\mod
4$. Assume that $(c,x+d)=1$ or $(d_0,x+c)=1$. If $4\mid xy$, then
$$(2a+\sqrt{4a^2+1})^{\f{p-1}4}\e\cases
(-1)^{\f d4+\f a2y+\f {xy}4}\mod p&\t{if $2\mid a$,}
\\(-1)^{\f d4+\f y4}\mod p&\t{if $2\nmid a$ and $4\mid y$,}
\\(-1)^{\f{a-1}2+\f x4}\f dc\mod p&\t{if $2\nmid a$ and $4\mid x$.}
\endcases$$
\endpro
Proof . This is immediate from (6.1), (6.2) and Theorem 7.1.
 \par\q
 \par From Theorem 7.1 and (6.3) we deduce the following result.
 \pro{Theorem 7.2} Let $a\in\Bbb Z$, and let $p\e 1\mod 8$ be a prime
such that $p=c^2+d^2=x^2+(4a^2+1)y^2$ with $c,d,x,y\in\Bbb Z$, $c\e
1\mod 4$ and $p\not=4a^2+1$. Suppose $d=2^rd_0$ with $d_0\e 1\mod
4$. Assume that $(c,x+d)=1$ or $(d_0,x+c)=1$. Then
$$p\mid U_{\f{p-1}8}(4a,-1)\iff \f{p-1}8\e\cases \f d4+\f
a2y+\f{xy}4\mod 2&\t{if $2\mid a$,}\\\f d4+\f y4\mod 2&\t{if $2\nmid
a$.}\endcases$$
\endpro

\par\q
\newline
{\bf Remark 7.1}  We conjecture that the condition $(c,x+d)=1$ or
$(d_0,x+c)=1$ in Theorems 6.1-6.4 and 7.1-7.2 can be canceled. See
also [S5, Conjectures 9.4, 9.10, 9.11 and 9.14].

\enddocument
\bye